\documentclass{amsart}
\usepackage{amssymb}
\usepackage{amscd}
\usepackage{graphics}

\usepackage[
hyperref,
nonotes
]{degt}

\def\ACC-{ACC-\penalty0\relax}

\def\ba{\bold a}
\def\bb{\bold b}
\def\bc{\bold c}

\def\be{\bold e}
\def\bv{\bold v}
\def\bB{\bold B}
\def\bL{\bold L}
\def\bP{\bold P}
\def\bS{\bold S}
\def\bT{\bold T}
\let\sset\bS
\def\tS{(\sset\oplus\Z h)\tilde\,}

\def\bx{\bar x}
\def\by{\bar y}
\def\bz{\bar z}

\def\barP{\bar P}
\def\barX{\bar X}
\let\t\tilde
\def\tP{\t P}
\def\tQ{\t Q}
\def\tF{\t F}
\def\tE{\t E}

\def\FF{\Bbb F_}

\let\poly\Pi
\let\cone\CalC
\let\pencil\CalP
\let\fiber=F

\def\height{\QOPNAME{ht}}
\def\genus{\QOPNAME{genus}}
\def\open{^\circ}

\let\Bertini\Gb

\let\Gm\mu
\let\Gn\nu

\def\DP{\GROUP{DP}}

\def\smax{\suptext{top}}
\def\smax{^\top\!}
\def\2{^{\prime2}}
\def\(#1){\{#1\}}

\let\polynomial\Cal
\def\pC{\polynomial C}

\def\pM{\polynomial M}
\def\pN{\polynomial N}
\def\pP{\polynomial P}
\def\pQ{\polynomial Q}
\def\pR{\polynomial R}
\def\pS{\polynomial S}
\def\pW{\polynomial W}

\let\sextic=D
\def\bsextic{\bar\sextic} 
\def\tsextic{\t\sextic}   
\def\psextic{\polynomial\sextic}
\def\bpsextic{\bar\psextic}

\let\rl=K
\def\brl{\bar\rl}

\def\prl{\polynomial\rl}
\def\bprl{\bar\prl}

\let\sect=L

\let\set\frak
\def\sQ{\set Q}

\def\black{$\bullet$}
\def\white{$\circ$}
\def\cross{$\scriptstyle\times$}

\newcounter{line}

\def\lineref#1{$\singset{#1}$, \autoref{l.#1}}

\let\PLUS+
\def\FRAC#1//#2>{\frac{#1}{#2}}
\def\LPAR{\left(}
\def\RPAR{\right)}
\def\EPS#1{\epsilon}
\let\leftbracket[
{\catcode`\*\active \catcode`\[\active \catcode`\+\active\catcode`\I\active
\catcode`\(\active\catcode`\)\active
\let\+\relax
\gdef\deMaple{%
 \catcode`\*\active \catcode`\[\active \catcode`\+\active\catcode`\I\active
 \catcode`\(\active\catcode`\)\active
 \def*{}\def[##1]{_{##1}}\let\[\leftbracket
 \def+{{}\PLUS{}}\let I=\epsilon\let\<\FRAC
 \let(\LPAR\let)\RPAR\def\{{\left.}\def\}{\right.}\let\e\EPS
 \def\next{\noalign{\vspace{-3pt}}&\phantom{{}:={}}}}%
\catcode`\_\active
\gdef\deZ{\deMaple\catcode`\_\active\def_Z{\epsilon}}%
}

\theoremstyle{definition}
 \addtheorem{Notation}{notation}

\title{On the Artal--Carmona--Cogolludo construction}

\author{Alex Degtyarev}

\address{%
Department of Mathematics\\
Bilkent University\\
06800 Ankara, Turkey}

\email{degt@fen.bilkent.edu.tr}



\begin{document}

\begin{abstract}
We derive explicit defining equations for
a number
of irreducible maximizing plane
sextics with double singular points only.
For most real curves, we also compute the fundamental group of the
complement; all groups found are abelian, which suffices to complete the
computation of the groups of all non-maximizing irreducible sextics.
As a by-product, examples of Zariski pairs in the strongest possible sense
are constructed.
\end{abstract}

\maketitle

\section{Introduction}

During
the last dozen of years, the geometry and topology of
singular complex plane
projective curves of degree six (\emph{plane sextics} in the sequel) has been
a subject of substantial interest.
Due to the fast development, there seems to be no good contemporary survey;
I can only suggest~\cite{degt:book} for a few selected topics and a number of
references.
Apart from the more subtle geometric properties that some special classes of
sextics may possess, the principal questions seem to be
\roster*
\item
the equisingular deformation classification of sextics,
\item
the fundamental group $\pi_1(\Cp2\sminus\sextic)$ of the complement,
\item
the defining equations.
\endroster
The last one seems more of a practical interest: the defining equations may
serve as a tool for attacking other problems.
However, the equations may also shed light on the arithmetical properties of
the so called \emph{maximizing} (\ie, those with the maximal total Milnor
number $\Gm=19$) sextics, as such curves are rigid (have discrete moduli
spaces) and are defined over algebraic number fields,
see~\cite{Persson:sextics}.

At present, the work is mostly close to its completion, at least for
\emph{irreducible} sextics. (Reducible sextics are too large in number on the
one hand and seem less interesting on the other.)
This paper bridges some of the remaining gaps.

In fact,\mnote{new paragraph; citations added}
the development of the subject is so fast that new results appear and become
available before old ones are published. Thus, the new
papers~\cite{degt:geography} and~\cite{Orevkov:equations} substantially
complement and complete the results of the present work. For the reader's
convenience, these new findings are either cited next to
or incorporated into the corresponding
statements.

All sextics with at least one triple or more complicated singular point,
including non-simple ones, are completely covered in~\cite{degt:book}
(the combinatorial approach used there seems more effective than the
defining equations); for this reason, such curves
are almost ignored
in this paper.
Thus, modulo a few quite
reasonable conjectures, which have mostly been proved, it
remains to study a few maximizing
sextics with double singular points only.

The classification of maximizing sextics is known:
it can be obtained
from~\cite{degt:JAG,Yang} (a reduction to an arithmetical problem),
\cite{Yang} (a list of the sets of singularities), and
\cite{Shimada:maximal} (the deformation classification).
The resulting list of \emph{irreducible}
maximizing sextics with double singular
points only is found in~\cite{degt:tetra}:
altogether, there are $39$ sets of
singularities realized by $42$ real and $20$ pairs of complex conjugate
curves.

Some of these sets of singularities have been studied and their defining
equations and fundamental groups are known, see~\cite{degt:tetra} for
references.
Here, we obtain the equations for $21$ set of singularities, leaving
only six sets unsettled, see~\eqref{eq.unknown}.
Then, we try to derive a few topological and
arithmetical consequences.

\subsection{Principal results}\label{s.results}
Strange as it seems, the main result of the paper does not appear in the
paper. We compute explicit defining equations for most irreducible maximizing
plane sextics with double singular points only.
Unfortunately, many equations are too complicated, and it seems neither
possible nor meaningful to reproduce them in a journal article. In both human
and machine readable form they can be downloaded from my web
page~\cite{degt:equations}; here, in \autoref{S.computation}, we outline the
details of the computation and provide information that is just sufficient to
recover the equations using the rather complicated formulas of
\autoref{S.Bertini}.
Note\mnote{new remark} that \cite{degt:equations} incorporates
as well the results
of~\cite{Orevkov:equations}, thus containing the defining equations of
\emph{all} irreducible maximizing sextics with double singular points only.

\table
\caption{Sextics considered in the paper}\label{tab.sextics}
\def\torus{(torus type)}
\def\spec#1{($\DG{#1}$-sextic)}
\def\torus{\spec6}
\def\spec#1{($\DG{#1}$)}
\def\*{\rlap{$^*$}}
\def\1{(1,0)}
\def\countr{0}
\def\countc{0}
\def\savecount#1#2{\count0=#1 \advance\count0 #2\xdef#1{\the\count0}}
\def\savecounts(#1,#2){\savecount\countr#1\savecount\countc#2$(#1,#2)$}
\def\FIRST#1.#2\\{\gdef\LABEL{#2}\space\strut\hss
 \refstepcounter{line}\theline.\label{l.\LABEL}}
\def\FIRST#1.#2\\{\setcounter{line}{#1}\addtocounter{line}{-1}%
 \gdef\LABEL{#2}\space\strut\hss
 \refstepcounter{line}\theline.\label{l.\LABEL}}
\hbox to\hsize{\hss\vbox{\obeylines\let
\cr%
\halign{\FIRST#\\%
 &\quad\singset{\LABEL}#\hss%
 &\quad\hss\expandafter\savecounts#\hss%
 &\quad\hss$#$\hss%
 &\quad#\hss\space
\noalign{\hrule\vspace{3pt}}%
\omit\strut\hss$\#$\ &\omit\quad Singularities\hss&\omit\quad\hss$(r,c)$\hss%
 &\pi_1&References, remarks
\noalign{\vspace{2pt}\hrule\vspace{2pt}}%
 4.A16+A3&&(2,0)&   \CG6&\eqref{eq.A16+A3}, see also \cite{Artal:trends}
 6.A15+A4&&(0,1)\*& \CG6&\eqref{eq.A15+A4}, see also \cite{Artal:trends}
 7.A14+A4+A1&&(0,3)&&\eqref{eq.A14+A4+A1}
10.A13+A6&&(0,2)&   &\eqref{eq.A13+A6}
11.A13+A4+A2&&(2,0)&\CG6&\eqref{eq.A13+A4+A2}
12.A12+A7&&(0,1)&   &\eqref{eq.A12+A7}
13.A12+A6+A1&&(1,1)&\CG6&\eqref{eq.A12+A6+A1}, see \autoref{rem.real.only}
14.A12+A4+A3&&\1&\CG6&\eqref{eq.A12+A4+A3}
16.A11+2A4&&(2,0)&\CG6&\eqref{eq.A11+2A4}, see \autoref{rem.one.only}
18.A10+A9&&(2,0)\*&\CG6&\eqref{eq.A10+A9}
19.A10+A8+A1&&(1,1)&\CG6&\eqref{eq.A10+A8+A1}, see \autoref{rem.real.only}
20.A10+A7+A2&&(2,0)&\CG6&\eqref{eq.A10+A7+A2}
21.A10+A6+A3&&(0,1)&&\eqref{eq.A10+A6+A3}
23.A10+A5+A4&&(2,0)&\CG6&\eqref{eq.A10+A5+A4}
24.A10+2A4+A1&&(1,1)&&\eqref{eq.A10+2A4+A1}, \eqref{eq.A10+2A4+A1'}
25.A10+A4+A3+A2&&\1&\CG6&\eqref{eq.A10+A4+A3+A2}
27.A9+A6+A4&&(1,1)\*&\CG6&\eqref{eq.A9+A6+A4}, see \autoref{rem.real.only}
30.A8+A7+A4&&(0,1)&&\eqref{eq.A8+A7+A4}
31.A8+A6+A4+A1&&(1,1)&\CG6&\eqref{eq.A8+A6+A4+A1}, see \autoref{rem.real.only}
34.A7+2A6&&(0,1)&&\eqref{eq.A7+2A6}
35.A7+A6+A4+A2&&(2,0)&\CG6&\eqref{eq.A7+A6+A4+A2}
\noalign{\vspace{2pt}\hrule\vspace{3pt}\footnotesize%
\hbox{\strut$^*$\ %
These sets of singularities are realized by reducible sextics as well}%
}\crcr}}\hss}
\endtable

We deal with
\emph{maximizing irreducible non-special sextics} (see the definition prior
to \autoref{th.special}). The sets of singularities for which equations
are obtained are listed in \autoref{tab.sextics}; for consistency, we retain
the numbering introduced in~\cite{degt:tetra}.
Also listed are the number of curves realizing each set of singularities (in
the form $(r,c)$, where $r$ is the number of real curves and $c$ is the
number of pairs of complex conjugate ones), the fundamental group
$\pi_1:=\pi_1(\Cp2\sminus\sextic)$, when known, and references to the
equations, other sources, and remarks.

The equations obtained are used to make a few observations stated in the four
theorems below; they concern the minimal fields of definition
(\autoref{th.minimal}), the fundamental group of the complement
(\autoref{th.pi1}), and a few examples of
the so-called arithmetic Zariski pairs
(Theorems~\ref{th.Zariski} and~\ref{th.Zariski.weak}).
It seems feasible that, with appropriate modifications, these statements
would extend to all irreducible maximizing sextics.

\theorem\label{th.minimal}
Let $n:=r+2c$ be the total number of irreducible sextics realizing
a
maximizing\mnote{statement edited incorporating \cite{Orevkov:equations}}
set of singularities~$\sset$
with double singular points only.
Then, the $n$
curves are
defined over an algebraic number
field~$\Bbbk$, $[\Bbbk:\Q]=n$\rom;
they differ by the $n$ embeddings $\Bbbk\into\C$.
If $n>2$, the Galois closure of~$\Bbbk$ has Galois group~$\DG{2n}$.
This field~$\Bbbk$ is
\emph{minimal} in the sense that
it is contained in the coefficient field of any defining polynomial.
\endtheorem

This
theorem is proved in \autoref{proof.minimal}, and the minimal fields of
definition are described in \autoref{S.computation}
together with the equations, see references in
\autoref{tab.sextics}. (Unless stated otherwise, $\Bbbk$ is the minimal field
containing the parameters listed in \autoref{S.computation}.)
This proof works as well for the few sextics with triple singular points
mentioned below, see~\eqref{eq.triple.point},
and, probably, for most other maximizing sextics.
In particular,\mnote{new remark} it works for the equation newly found
in~\cite{Orevkov:equations} (see~\cite{degt:equations} for details),
and this fact is incorporated into the statement.
Another proof, using
the concept of
\emph{dessins d'enfants}, is
discussed in
\autoref{proof.alt}; it leads to somewhat disappointing consequences,
see \autoref{rem.bad.dessin}.

\theorem\label{th.pi1}
With two exceptions,\mnote{statement edited incorporating
\cite{Orevkov:equations}}
the fundamental group of a \emph{real} maximizing non-special sextic with
double singular points only is~$\CG6$.
The exceptions are
the real curve realizing \lineref{A10+2A4+A1}
and one of the two curves realizing \lineref{A11+2A4},
see \autoref{rem.one.only} below.
\endtheorem

In the two exceptional cases, the fundamental group is unknown
(at least, to the author).
I expect that these groups are also abelian, as well as those of the non-real
curves in \autoref{tab.sextics}.
The\mnote{extended}
proof of this theorem is partially based upon~\cite{Orevkov:equations};
it is explained in \autoref{proof.pi1}, and all technical details are found
in~\cite{degt:equations}.

\remark\label{rem.real.only}
The sets of singularities
\lineref{A12+A6+A1},
\lineref{A10+A8+A1},
\lineref{A9+A6+A4},
and \lineref{A8+A6+A4+A1}
are realized by three Galois conjugate curves each. In
each case, only one of the three curves is real, and
\emph{only for this real curve}
the fundamental group $\pi_1=\CG6$ has been computed.
\endremark

\remark\label{rem.one.only}
The set of singularities
\lineref{A11+2A4}
is realized by two Galois conjugate curves. Both curves
are real, but the group $\pi_1=\CG6$ is computed for one of them
only; for the other curve, the presentation obtained is incomplete and I
cannot assert that the group is finite.
\endremark

The following corollary of \autoref{th.pi1} relies on the deformation
classification of irreducible sextics, which is now completed,
see~\cite{degt:geography};\mnote{citation added}
the proof
will appear elsewhere.

\corollary[see \cite{degt:geography}]
Let $\sextic\subset\Cp2$ be a non-maximizing non-special irreducible simple
plane sextic. Then,
unless the set of singularities of~$\sextic$ is
\[*
\singset{2D7+2A2},\quad
\singset{D7+D4+3A2},\quad
\singset{2D4+4A2},\quad\text{or}\quad
\singset{2A4+2A3+2A2},
\]
one has $\pi_1(\Cp2\sminus\sextic)=\CG6$.
\pni
\endcorollary

In the four exceptional cases, the
groups are known to be
non-abelian: they are $\SL(2,\FF5)\odot\CG{12}$
for the last curve and
$\SL(2,\FF3)\times\CG2$ for the three others.
Here,\mnote{explained}
the notation $\SL(2,\FF5)\odot\CG{12}$ stands for the \emph{central product},
\ie, the direct product $\SL(2,\FF5)\times\CG{12}$ with the center
$\CG2\subset\SL(2,\FF5)$ identified with $\CG2\subset\CG{12}$.

Although special care has been taken to avoid sextics with triple singular
points,
see \autoref{rem.forbidden},
some of them do appear in the computation. These are the curves
realizing the following eight sets of singularities:
\[
\def\entry#1,{\strut\singset{#1}, see~\eqref{eq.#1}}
\def\ha#1{\vcenter{\halign{\entry##\hss\cr#1\crcr}}}
\hbox{$\ha{%
E6+A13,,\cr
E6+A10+A3,,\cr
E6+A7+A6,,\cr
D9+A10,,\cr
}\qquad\qquad\ha{%
D9+A6+A4,,\cr
D5+A14,,\cr
D5+A10+A4,,\cr
D5+A8+A6,.\cr
}$}
\label{eq.triple.point}
\]
Their equations are also described in \autoref{S.computation}, and
the conclusion of \autoref{th.minimal} extends to these curves literally,
together with the proof.
The fundamental groups of all these curves are abelian, see~\cite{degt:book}.

In this paper, we confine ourselves to the sextics that can be obtained from
a pencil of cubics with at most four basepoints (see \autoref{s.idea} for the
explanation).
The remaining six sets of singularities are\mnote{edited}
\[
\def\entry#1.#2,{\strut#1.\ \singset{#2}}
\def\ha#1{\vcenter{\halign{\entry##\hss\cr#1\crcr}}}
\hbox{$\ha{%
15.A12+A4+A2+A1,,\cr
22.A10+A6+A2+A1,,\cr
26.A10+A4+2A2+A1,,\cr
}\qquad\qquad\ha{%
36.A7+2A4+2A2,,\cr
38.2A6+A4+A2+A1,,\cr
39.A6+A5+2A4,.\cr
}$}
\label{eq.unknown}
\]
For these curves,
the models used in \autoref{S.computation} are not
applicable: pencils
with more than four basepoints need to be considered and the
computation seems to become much more involved.
In~\cite{degt:equations}, the defining equations of these six curves are
derived from the parametric equations found
in~\cite{Orevkov:equations}.\mnote{extended}

We conclude with a few examples of the so-called \emph{Zariski pairs},
see~\cite{Artal:Zariski}.
Roughly, two plane curves $D_1,D_2\subset\Cp2$
constitute a Zariski pair if they are combinatorially equivalent but
topologically distinct, in the sense that may vary from problem to problem.
In \autoref{th.Zariski}, we use the strongest combinatorial equivalence
relation (the curves are Galois conjugate;
in the terminology of~\cite{Shimada:Zariski},
the Zariski pairs are \emph{arithmetic}) and the weakest topological one
(the complements $\Cp2\sminus D_i$,\mnote{`()' removed} $i=1,2$ are not properly
homotopy equivalent).
In \autoref{th.Zariski.weak}, the topological relation is
slightly stronger.
The first example of Galois conjugate
but not homeomorphic algebraic varieties is
due to Serre~\cite{Serre:Zariski}.
Still, very few other examples are known;
a brief survey of the subject, including arithmetic Zariski pairs on plane
curves,
is contained in~\cite{Shimada:Zariski}.
For Zariski pairs in general, see~\cite{Artal:survey}.

\theorem\label{th.Zariski}
Let~$\sset$ be one of the following twelve sets of singularities\rom:
\let\lref\lineref
\roster
\item\label{Zariski.1}
\singset{D5+A10+A4},
\singset{A18+A1},
\singset{A16+A2+A1},
\lref{A12+A6+A1},
\singset{A12+A4+A2+A1},
\lref{A10+A8+A1},
\singset{A10+A6+A2+A1},
\lref{A10+2A4+A1},
\lref{A8+A6+A4+A1}\rom;
\item\label{Zariski.2}
\singset{E6+A10+A3},
\lref{A16+A3},
\lref{A10+A9},
\lref{A10+A7+A2},
\singset{A10+A4+2A2+A1},
\lref{A7+A6+A4+A2}.
\endroster
In case~\ref{Zariski.1}, with $(r,c)=(1,1)$,
let~$\sextic_1$, $\sextic_2$ be
a real and a
non-real sextic realizing~$\sset$\rom;
in case~\ref{Zariski.2}, with $(r,c)=(2,0)$, let $\sextic_1$, $\sextic_2$
be the two real sextics.
Then $(\sextic_1,\sextic_2)$ is a Zariski pair in the following strongest
sense\rom:
the two curves
are Galois conjugate, but the complements
$\Cp2\sminus\sextic_i$, $i=1,2$, are not properly homotopy equivalent.
\endtheorem

For the four sets of singularities
\singset{A18+A1},
\lineref{A16+A3},
\singset{A16+A2+A1}, and
\lineref{A10+A9}, the fact that the spaces $\Cp2\sminus D_i$, $i=1,2$, are
not homeomorphic was
originally established in~\cite{Shimada:Zariski}, and for
\autoref{th.Zariski} we use essentially the same topological invariant.

\theorem\label{th.Zariski.weak}
Let~$\sset$ be either \lineref{A13+A4+A2}, or \lineref{A10+A5+A4},
or \singset{A6+A5+2A4}.
Then the pair
$(\sextic_1,\sextic_2)$ of two distinct real sextics realizing~$\sset$ is a
Zariski pair in the following sense\rom: the two curves are Galois
conjugate, but the topological pairs
$(\Cp2,\sextic_i)$, $i=1,2$, are not homotopy equivalent.
In fact,\mnote{moved} a slightly stronger statement holds\rom: the pairs
\[
(\Cp2\sminus\Sing\sextic_i,\sextic_i\sminus\Sing\sextic_i),\quad
i=1,2
\label{eq.pairs}
\]
are not properly homotopy equivalent,
where $\Sing\sextic_i$ stands for the set of all singular points of~$\sextic_i$.
\endtheorem

These
theorems are proved in \autoref{proof.Zariski} and
\autoref{proof.Zariski.weak}.
For the sets of singularities as in \autoref{th.Zariski}\iref{Zariski.2}
and \autoref{th.Zariski.weak},
we can also state that the Zariski pairs are
\emph{$\pi_1$-equivalent}, \ie, the fundamental groups
$\pi_1(\Cp2\sminus\sextic_i)$, $i=1,2$, are isomorphic
(as they are both~$\CG6$).
In fact, the spaces $\Cp2\sminus\sextic_i$ are homotopy equivalent,
see \autoref{prop.homotopy}, but they are not homeomorphic!
Probably, this conclusion
(the homotopy equivalence of the complements)
also holds for \autoref{th.Zariski}\iref{Zariski.1}.\mnote{next paragraph
removed, incorporated into the statements}


\subsection{The idea}\label{s.idea}
The main tool used in the paper is the
Artal--Carmona--Cogolludo construction
(\emph{\ACC-construction} in the sequel)
developed in~\cite{Artal:trends}.
This construction is outlined in \autoref{s.ACC}.
We confine ourselves to generic (in the sense of
\autoref{def.generic}) irreducible non-special sextics with double singular
points only, see \autoref{conv.sextic};
using the theory of $K3$-surfaces,
we show that the ramification locus of the \ACC-model
of such a curve is
irreducible and with $\bA$ type singularities only, see \autoref{s.D}.
Furthermore, we describe the singular fibers of the corresponding
Jacobian elliptic
surface~$Y$, see \autoref{s.fibers}, and prove that,
geometrically, the blow-down
map $Y\to\Cp2$
is as shown in \autoref{fig.blow-down} on \autopageref{fig.blow-down}.
This fact eliminates the need in the tedious case-by-case analysis of the
possible configurations of divisors (\cf.~\cite{Artal:trends}), and the
construction of the \ACC-models of all
maximizing sextics
becomes a relatively easy task.

The other key ingredient is Moody's paper~\cite{Moody}.
With a little effort (and \Maple) its results can be extended to explicit
formulas for the rational
two-to-one map $\Cp2\dashrightarrow\Sigma_2$ related to
the \ACC-model, see \autoref{s.map}; they are used to pass from the
defining equations of the \ACC-models in~$\Sigma_2$ to those of the original
sextics.
(In fact, we incorporate Moody's formulas from the very beginning and
describe the \ACC-models in terms of pencils of cubics; it may be due to this
fact that, in the six missing cases,
the equations on the parameters are too complicated
to be solved or even to be written down.)

The fundamental groups are computed as suggested in~\cite{degt:tetra},
representing sextics as tetragonal curves. Alternatively, one could try to
use the \ACC-models, which are \emph{trigonal} curves and look simpler;
unfortunately, one would have to keep track of too many (three to four) extra
sections, which makes this approach about as difficult as the direct
computation, especially when the curve is not real.
Certainly, given equations, one can also use the modern
technology and compute the monodromy by brute force; however, at this stage I
prefer to refrain from a computer aided solution to a problem that is not
discrete in its nature.

The other theorems are proved in \autoref{S.proofs} by constructing
appropriate invariants.

\subsection{Contents of the paper}
In \autoref{S.K3}, after a brief introduction to the basic concepts
related to plane sextics,
we use the theory of $K3$-surfaces to describe the rational curves splitting
in the double covering ramified at a generic non-special irreducible sextic.
These results are used in \autoref{S.ACC}, where we introduce the \ACC-model
and show that the models of non-special
irreducible sextics are particularly simple.
In \autoref{S.Bertini}, we recall and extend the results of~\cite{Moody}
concerning the Bertini involution $\Cp2\dashrightarrow\Cp2$
and explain how these results apply to the \ACC-construction.
In \autoref{S.computation}, the details of deriving the defining equations of
maximizing sextics are outlined and their minimal fields of definition are
described.
Finally, in \autoref{S.proofs} we give formal proofs of
Theorems \ref{th.minimal}, \ref{th.pi1}, \ref{th.Zariski},
and~\ref{th.Zariski.weak} and make a few concluding remarks.
As a digression, we discuss the homotopy type of the complement of an
irreducible plane curve with abelian fundamental group, see
\autoref{s.homotopy}.

\subsection{Acknowledgements}
I am grateful to Igor Dolgachev, who kindly explained to me
an alternative approach to the treatment of the Bertini involution, to
Alexander Klyachko, who patiently introduced me to the
more practical aspects of
Galois theory, to Sergey Finashin, who brought to my attention
paper~\cite{Dyer.Sieradski},
and to Stepan Orevkov, who generously shared his
results~\cite{Orevkov:equations}.\mnote{new acknowledgement}

\section{The covering $K3$-surface}\label{S.K3}

The principal goal of this section is \autoref{th.rational},
which describes
the rational curves in the $K3$-surface ramified at a generic
(see \autoref{def.generic} below)\mnote{ref added}
non-special irreducible sextic.

\subsection{Terminology and notation}
A plane sextic $\sextic\subset\Cp2$ is called \emph{simple}
if all its singular points are
simple, \ie, those of type $\bA$--$\bD$--$\bE$, see~\cite{Durfee}.
Given a sextic~$\sextic$, we denote by $P_i\in\Cp2$ its singular points.

We will also use the classical concept of infinitely near points: given a
point~$P$ in a surface~$S$, all points~$Q$ in the exceptional divisor in the
blow-up $S(P)$ of~$S$ at~$P$ are said to be \emph{infinitely near} to~$P$
(notation $Q\to P$).
A curve $D\subset S$ \emph{passes} through a point $Q\to P$
(notation $Q\in D$)
if $D$ passes
through~$P$ and
the strict transform of~$D$ in $S(P)$ passes through~$Q$.
Similarly, a
point $Q\to P$ is \emph{singular} for~$D$ if it is singular for
the strict transform of~$D$ in $S(P)$.
This construction can be iterated and one can consider sequences
$\ldots\to Q'\to Q\to P$ of infinitely near points.
Starting from level~$0$ for the points of the original plane~$\Cp2$,
we define the \emph{level} of an infinitely near point \via\
$\QOPNAME{level}(Q)=\QOPNAME{level}(P)+1$ whenever $Q\to P$.

For a sextic $\sextic\subset\Cp2$ with $\bA$-type singular points only, denote
by $\DP(\sextic)$
the set of all \emph{double points} of~$\sextic$, including infinitely near.
This set is a union of disjoint maximal (with respect to inclusion) chains,
each
singular point~$P_i$ of type~$\bA_p$ giving rise to a chain
$Q_r\to\ldots\to Q_1=P_i$ of length $r=[\frac12(p+1)]$.
A subset $\sQ\subset\DP(\sextic)$ is called \emph{complete} if, whenever $Q\in\sQ$
and $Q\to Q'$, also $Q'\in\sQ$.

\subsection{The homological type}\label{s.K3}
Given a simple
sextic~$\sextic$, denote by $X:=X_\sextic$ the minimal resolution of singularities of
the double covering of the plane~$\Cp2$ ramified at~$\sextic$. It is a
$K3$-surface, see, \eg,~\cite{Persson:sextics}.
With a certain abuse of the language, $X$ is referred to as the
\emph{covering $K3$-surface}.

Denote by $\bL:=H_2(X)\cong2\bE_8\oplus3\bU$ the intersection index lattice
of~$X$ (where $\bU=\Z\bold u_1+\Z\bold u_2$, $\bold u_1^2=\bold u_2^2=0$,
$\bold u_1\cdot\bold u_2=1$ is the \emph{hyperbolic plane}).
Let $h\in\bL$, $h^2=2$, be the class of a hyperplane section (pull-back of a
line in~$\Cp2$), and let $\bP_i\subset\bL$ be the lattice spanned by the
classes of the exceptional divisors over a singular point~$P_i$ of~$\sextic$.
It is a negative definite even root lattice of the same name
$\bA$--$\bD$--$\bE$ as the type of~$P_i$, and $\rank\bP_i=\Gm(P_i)$ is the
Milnor number.
This lattice has a \emph{canonical basis}, constituted by the classes of the
exceptional divisors. The basis vectors are the walls of a single Weyl
chamber; they can be identified with the vertices of the Dynkin graph
of~$\bP_i$.

We will consider the sublattice $\sset:=\bigoplus_i\bP_i$, referred to as the
\emph{set of singularities} of~$\sextic$,
the obviously
orthogonal sum $\sset\oplus\Z h$, and its \emph{primitive hull}
\[*
\tS:=\bigl((\sset\oplus\Z h)\otimes\Q\bigr)\cap\bL.
\]
The sequence of lattice extensions
\[
\sset\subset\sset\oplus\Z h\subset\bL
\label{eq.h-type}
\]
is called the \emph{homological type} of~$\sextic$. Clearly, $\sset$ is an
even\mnote{`even'}
negative definite lattice and $\rank\bS=\Gm(\sextic)$ is the total
Milnor number of~$\sextic$. Since $\Gs_-\bL=19$, one has $\Gm(\sextic)\le19$,
see~\cite{Persson:sextics}. A simple sextic~$\sextic$ with $\Gm(\sextic)=19$
is called \emph{maximizing}. Note that both the inequality and the term apply
to simple sextics only.

An irreducible sextic $\sextic\subset\Cp2$ is called \emph{special},
or \emph{$\DG{2n}$-special},
see~\cite{degt:Oka}, if
its fundamental group $\pi_1(\Cp2\sminus \sextic)$ admits a dihedral
quotient $\DG{2n}$, $n\ge3$. If this is the case,
one has $n=3$, $5$, or~$7$, and $\DG6$-special sextics are those of
\emph{torus type},
see~\cite{degt:Oka}.
By definition, the fundamental groups of special sextics are not abelian.

\theorem[see~\cite{degt:Oka}]\label{th.special}
A simple sextic~$\sextic$ is irreducible and non-special if and only if
$\sset\oplus\Z h\subset\bL$ is a primitive sublattice, \ie,
$\sset\oplus\Z h=\tS$.
\pni
\endtheorem

Since both~$\sset$ and $\Z h$ are generated by algebraic curves (and since
the N\'{e}ron--Severi lattice $\NS(X)$ is primitive in~$\bL$), one has
$\tS\subset\NS(X)$.

\definition\label{def.generic}
A simple sextic $\sextic\subset\Cp2$ is called \emph{generic}
if $\tS=\NS(X)$.
\enddefinition

In each equisingular stratum of the space of simple sextics,
generic ones form a dense Zariski open subset.
A maximizing sextic is always generic.

Let $\tau\:X\to X$ be the deck translation of the ramified covering
$X\to\Cp2$. This automorphism induces an involutive autoisometry
$\tau_*\:\bL\to\bL$.

\lemma\label{lem.tau*}
The induced autoisometry $\tau_*\:\bL\to\bL$ acts as follows\rom:
$\tau_*(h)=h$\rom; the restriction of~$\tau_*$ to~$\bP_i$ is
induced by the symmetry~$s_i$ of the Dynkin graph of~$\bP_i$
\rom(in the canonical basis\rom), where
\roster*
\item
$s_i$ is the only nontrivial symmetry
if $\bP_i=\bA_{p\ge2}$, $\bD\subtext{odd}$, or~$\bE_6$, and
\item
$s_i$ is the identity otherwise\rom;
\endroster
the restriction of~$\tau_*$ to $(\sset\oplus\Z h)^\perp$ is $-\id$.
\endlemma

\proof
The action of~$\tau_*$ on $\sset\oplus\Z h$ is given by a simple computation
using the minimal resolution of the singularities of~$\sextic$ in $\Cp2$.
For the last statement,
since $X/\tau=\Cp2$ is rational, $\tau$ is anti-symplectic, \ie,
$\tau_*(\Go)=-\Go$ for the class~$\Go$ of a holomorphic form on~$X$.
Since also $\tau_*$ is defined over~$\Z$,
the $(-1)$-eigenspace of $\tau_*$ contains the minimal \emph{rational}
subspace $V\subset\bL\otimes\Q$ such that $\Go\in V\otimes\C$. On the other
hand, $\tau_*$ is invariant under equisingular deformations and,
deforming~$\sextic$ to a generic sextic, one has
$V=(\sset\oplus\Z h)^\perp\otimes\Q$.
(Recall that $\NS(X)=\Go^\perp\cap\bL$.)
\endproof

\subsection{Rational curves in~$X$}\label{s.rational}
The goal of this section is the following theorem, which is proved at the end
of the section.

\theorem\label{th.rational}
Let $\sextic\subset\Cp2$ be a
generic irreducible non-special sextic with $\bA$-type
singularities only, and let $X$ be the covering $K3$-surface.
Let, further, $R\subset X$ be a nonsingular rational curve whose projection
$\bar R\subset\Cp2$ is a curve of degree at most~$3$.
Then the projection $R\to\bar R$ is two-to-one
\rom(in other words, $R$ is the pull-back
of the strict transform of~$\bar R$\rom), and $\bar R$ is one of the
following\rom:
\roster
\item\label{rational.1}
a line through a complete pair $\sQ_2\subset\DP(\sextic)$\rom;
\item\label{rational.2}
a conic through a complete quintuple~$\sQ_5\subset\DP(\sextic)$\rom;
\item\label{rational.3}
a cubic through a complete septuple~$\sQ_7\subset\DP(\sextic)$ with a
double point at a distinguished point $P\in\sQ_7$ of level zero.
\endroster
Conversely, given a complete set~$\sQ_2$, $\sQ_5$ or pair $P\in\sQ_7$ as
above, there is a unique, respectively, line, conic, or cubic~$\bar R$ as in
items~\ref{rational.1}--\ref{rational.3}. This curve~$\bar R$ is
irreducible, and the pull-back of its strict transform is a nonsingular
rational curve in~$X$.
\endtheorem

\corollary\label{cor.generic}
Under the hypotheses of \autoref{th.rational}, the configuration
$\DP(\sextic)$ is \emph{almost del Pezzo}, in the sense that
\roster
\item\label{delPezzo.1}
there is no line passing through three points\rom;
\item\label{delPezzo.2}
there is no conic passing through six points\rom;
\item\label{delPezzo.3}
there is no cubic passing through eight points and singular at one of them.
\endroster
More generally, there is no line or conic whose local intersection index
with~$\sextic$ at each intersection point is even.
\done
\endcorollary

\remark
If $\sextic$ is a special sextic, there \emph{are} conics (not necessarily
irreducible) passing through some
complete sextuples $\sQ_6\subset\DP(\sextic)$,
see~\cite{degt:Oka}. Thus, the existence of such conics is yet another
characterization of \emph{generic} irreducible special sextics.
See~\cite{Shimada:splitting} for more details.
\endremark

We precede the proof of \autoref{th.rational} with a few observations.

Let $\be_1,\ldots,\be_k$ be the canonical basis for a summand of~$\sset$ of
type~$\bA_k$.
For an element $a\in\bA_k$, $a=\sum a_i\be_i$,
let $a_0=a_{k+1}=0$ and denote
$d_i=a_{i}-a_{i-1}$, $i=1,\ldots,k+1$.
Recall that $\bA_k$ is the orthogonal complement
$(\bv_1+\ldots+\bv_{k+1})^\perp$ in the lattice
$\bB_{k+1}:=\bigoplus_{i=1}^{k+1}\Z\bv_i$, $\bv_i^2=-1$,
so that $\be_i=\bv_{i}-\bv_{i+1}$, $i=1,\ldots,k$. In this notation, one has
$a=\sum d_i\bv_i$. Hence, $a^2=-\sum d_i^2$ and
$a\cdot a'=-\sum d_id_i'$ for another element $a'=\sum d_i'\bv_i$.
The following statement is straightforward.

\lemma\label{lem.Ak}
An element $a=\sum d_i\bv_i\in\bB_{k+1}$ is in~$\bA_k$ if and only if
$\sum d_i=0$.
Furthermore, $a\cdot\be_i\ge0$ for all $i=1,\ldots,k$ if and only if
$d_1\le d_2\le\ldots\le d_{k+1}$.
\done
\endlemma

\corollary\label{cor.A10}
The elements $a\in\bA_k$ such that $a^2\ge-10$ and $a\cdot\be_i\ge0$ for all
$i=1,\ldots,k$ are as follows\rom:
\roster*
\item
$\ba^q:=-\bv_1-\ldots-\bv_q+\bv_{k+2-q}+\ldots+\bv_{k+1}$
\rom($1\le q\le5$\rom)\rom:
$(\ba^q)^2=-2q$\rom;
\item
$\bb^q:=\ba^1+\ba^q$ \rom($1\le q\le 2$\rom)\rom: $(\bb^q)^2=-6-2q$\rom;
\item
$\bc^+:=-2\bv_1+\bv_k+\bv_{k+1}$ or $\bc^-:=-\bv_1-\bv_2+2\bv_{k+1}$\rom:
$(\bc^\pm)^2=-6$.
\done
\endroster
\endcorollary

\corollary\label{cor.A20}
The
\emph{symmetric} \rom(with respect to the
only nontrivial symmetry of the Dynkin
graph\rom)
elements $a\in\bA_k$ such that $a^2\ge-20$ and $a\cdot\be_i\ge0$ for all
$i=1,\ldots,k$ are
$\ba^q$ \rom($1\le q\le10$\rom),
$\bb^q$ \rom($1\le q\le7$\rom),
$\ba^2+\ba^q$ \rom($2\le q\le4$\rom),
and $2\ba^1+\ba^q$ \rom($1\le q\le2$\rom),
see \autoref{cor.A10} for the notation.
\done
\endcorollary

\proof[Proof of \autoref{th.rational}]
Recall a description of rational curves on a $K3$-surface~$X$,
see~\cite{Pjatecki-Shapiro.Shafarevich} or~\cite[Theorem~6.9.1]{DIK}.
Let $\cone:=\{x\in\NS(X)\otimes\R\,|\,x^2>0\}$ be the positive cone and
$\PP(\cone):=\cone/\R^*$
its projectivization. Consider the group~$G$ of motions of the
hyperbolic space $\PP(\cone)$ generated by the reflections against
hyperplanes orthogonal to vectors $v\in\NS(X)$ of square~$(-2)$ and let
$\poly\subset\PP(\cone)$ be the fundamental polyhedron of~$G$ containing the
class of a K\"{a}hler form $\rho\in\NS(X)\otimes\R$.
Denote by $\Delta_+(X)$ the set of vectors $v\in\NS(X)$ such that $v^2=-2$,
$v\cdot\rho>0$, and $v$ is orthogonal to a face of~$\poly$.
Then
{\em $\Delta_+(X)\subset\bL$
is precisely the set of homology classes realized by
nonsingular rational curves on~$X$.}
Each class $v\in\Delta_+(X)$ is realized by a unique such curve.

The set $\Delta_+:=\Delta_+(X)$
can be found step by step,
using Vinberg's algorithm~\cite{Vinberg:polyhedron} and taking for~$\rho$ a
small perturbation of~$h$. At Step~0, one adds to $\Delta_+$ the classes of
all exceptional divisors, \ie, the canonical basis for~$\sset$. Then, at
each step~$s$, $s>0$, one adds to~$\Delta_+$ all vectors $v\in\NS(X)$ such
that $v^2=-2$, $v\cdot h=s$ (there are finitely many such vectors), and
$v\cdot u\ge0$ for any $u\in\Delta_+$ with $u\cdot h<s$, \ie,
{\em $v$ has non-negative intersection with any vector added to $\Delta_+$ at
the previous steps.}
\latin{A priori}, $\Delta_+(X)$ may be infinite and this
process does not need to terminate; for further details and the termination
condition, see~\cite{Vinberg:polyhedron}.

Under the hypotheses of the theorem, $\NS(X)=\sset\oplus\Z h$ and all
vectors~$v$ used in Vinberg's algorithm are of the form $v=a+rh$, $r>0$,
$a\in\sset$, $a^2=-2r^2-2$.
In particular, all odd steps are vacuous and
the condition $v\cdot u\ge0$ for all $u\in\Delta_+$, $u\cdot h=0$ is
equivalent to the requirement that each component $a_i\in\bP_i$
of~$a$ should be as in \autoref{lem.Ak}.
With one exception (the curve~$\sextic$ itself, see below),
the rational curve~$R$ represented by $v=a+rh$ projects to a curve
$\bar R\subset\Cp2$ of degree~$2r$ or~$r$, depending on whether $v$ is or is
not $\tau_*$-invariant; the projection is two-to-one in the former case and
one-to-one in the latter. (Thus, the projection is two-to-one whenever
$\deg\bar R$ is odd, and only the case of conics needs special attention.)

Summarizing, we need to consider Steps~2, 4, and~6 of
Vinberg's algorithm, in the
last case confining ourselves to $\tau_*$-invariant vectors only.

Introduce the notation $\ba_i^q\in\bP_i$ \etc\. for the elements $\ba^q$
\etc\. as in
\autoref{cor.A10} in the lattice~$\bP_i$. Throughout
the rest of the proof, we assume
the convention
that distinct subscripts represent distinct indices;
for example, an expression
$\ba_i^1+\ba_j^1$ implies implicitly that $i\ne j$.

Added at Step~2 are all elements of the form $\ba_i^2+h$ (whenever
$\Gm(P_i)\ge3$) and $\ba_i^1+\ba_j^1+h$. These elements are in a
one-to-one correspondence with
complete pairs $\sQ_2\subset\DP(\sextic)$, and the
corresponding rational curves are the pull-backs of the lines as in
\autoref{th.rational}\iref{rational.1}.
In particular, we conclude that there are no conics
in $\Cp2$ whose pull-backs split into pairs of rational
curves in~$X$.

Added at Step~4 are all elements of the form $\sum\ba_{i_\Ga}^{q_\Ga}+2h$,
$\sum q_\Ga=5$.
These elements are in a one-to-one correspondence with complete quintuples
$\sQ_5\subset\DP(\sextic)$,
and the rational curves are the pull-backs of the conics as in
\autoref{th.rational}\iref{rational.2}.

Note that elements containing $\bb^q$ or $\bc^\pm$,
see \autoref{cor.A10}, cannot be added at Step~4.
Indeed, any such element would be one of the following:
\roster*
\item
$a=\bb_i^2+2h$; then $a\cdot(\ba_i^2+h)=-2<0$;
\item
$a=\bb_i^1+\ba_j^1+2h$; then $a\cdot(\ba_i^1+\ba_j^1+h)=-2<0$;
\item
$a=\bc_i^\pm+\ba_j^2+2h$; then $a\cdot(\ba_i^1+\ba_j^1+h)=-1<0$;
\item
$a=\bc_i^\pm+\ba_j^1+\ba_k^1+2h$; then $a\cdot(\ba_i^1+\ba_j^1+h)=-1<0$.
\endroster

Finally, we characterize the \emph{$\tau_*$-invariant} elements added at
Step~6.

An element of the form $\sum\ba_{i_\Ga}^{q_\Ga}+3h$, $\sum q_\Ga=10$, exists
(and then is unique) if and only if $\sextic$ has ten double points. In this case,
$\genus(\sextic)=0$ and the above element represents~$\sextic$ itself.

Elements of the form $\bb_i^q+\sum\ba_{i_\Ga}^{q_\Ga}+3h$, $q+\sum q_\Ga=7$,
\emph{are} added at Step~6.
(It is
immediate
that any such element has non-negative intersection
with all those added at Steps~2 and~4.)
Such elements are parametrized by pairs $P_i\in\sQ_7$, where
$\sQ_7\subset\DP(\sextic)$ is a
complete septuple and $P_i$ is a distinguished point of
level zero.
The corresponding rational curve is the pull-back of a cubic as in
\autoref{th.rational}\iref{rational.3}.
No other $\tau_*$-invariant element can be added at this step. Indeed,
by \autoref{cor.A20}, such an element would be one of the following:
\roster*
\item
$a=(\ba_i^2+\ba_i^q)+\ldots+3h$ ($2\le q\le4$); then $a\cdot(\ba_i^2+h)=-2<0$;
\item
$a=3\ba_i^1+\ba_j^1+3h$; then $a\cdot(\ba_i^1+\ba_j^1+h)=-2<0$;
\item
$a=(2\ba_i^1+\ba_i^2)+3h$; then $a\cdot(\ba_i^2+h)=-2<0$.
\endroster
This observation completes the proof of \autoref{th.rational}.
\endproof

\section{The \ACC-construction}\label{S.ACC}

In this section, we recall the principal results of~\cite{Artal:trends}
concerning the properties of the \ACC-construction
and describe the singular fibers and the ramification locus of the \ACC-model
of a generic non-special irreducible sextic.

\subsection{The construction\noaux{ (see~\cite{Artal:trends})}}\label{s.ACC}
Consider a sextic $\sextic\subset\Cp2$ with $\bA$-type singular points only and
fix a complete octuple
$\sQ_8=\{Q_1,\ldots,Q_8\}\subset\DP(\sextic)$. (Thus, we assume that $\sextic$ has
at least eight double points. If $\sextic$ is irreducible, this
assumption is equivalent
to the requirement that $\genus(\sextic)\le2$.)

Let $\pencil:=\pencil(\sQ_8)$ be the closure of the set
of cubics passing through all points
of~$\sQ_8$. As shown in~\cite{Artal:trends},
$\pencil$ is a pencil and
a generic member of~$\pencil$
is a nonsingular cubic; hence, $\pencil$ has nine basepoints: the points
$Q_1,\ldots,Q_8$
of~$\sQ_8$ and another
\emph{implicit}
point $Q_0$, which \latin{a priori} may be infinitely
close to some of $Q_i\in\sQ_8$.
Let $\sQ_8^*:=\sQ_8\cup\{Q_0\}$.
It follows that the result $Y:=\Cp2(\sQ^*_8)$ of the blow-up of the nine points
$Q_0,\ldots,Q _8$ is a
relatively minimal
rational Jacobian elliptic surface, the distinguished
section being the exceptional divisor~$\tQ_0$ over~$Q_0$.
With $\sQ_8$ (and hence $Q_0$ and~$Y$) understood, we use the notation $\t A$
for the strict transform in~$Y$ of a curve $A\subset\Cp2$.
Let also $\tQ_i\subset Y$ be the strict transform of the exceptional
divisor obtained by blowing up the basepoint~$Q_i$, $i=0,\ldots,8$, and let
$\t\pencil$ be the resulting elliptic pencil on~$Y$.

The fiberwise
multiplication by~$(-1)$ is an involutive automorphism
$\Bertini\:Y\to Y$, and the
quotient $Y/\Bertini$ blows down to the Hirzebruch surface~$\Sigma_2$,
\ie, geometrically ruled rational surface with an exceptional section~$E$ of
self-intersection~$(-2)$ (the image of~$\tQ_0$).
Conversely, $Y$ is recovered as the minimal resolution of singularities
of the double covering of~$\Sigma_2$ ramified at the exceptional section~$E$
and a certain \emph{proper trigonal curve}~$\brl$ (\ie, a reduced
curve disjoint from~$E$ and intersecting each fiber of the ruling at three
points).
This representation of the Jacobian elliptic surface~$Y$ is often referred to
as its \emph{Weierstra{\ss} model}.

\theorem[see~\cite{Artal:trends}]\label{th.Artal}
One has $\tsextic\cdot\tQ_0=0$ and $\tsextic\cdot\tF=2$, where $\tF$ is a generic
fiber of~$\t\pencil$.
Furthermore, one has $\Bertini(\tsextic)=\tsextic$\rom;
hence, the image
$\bsextic\subset\Sigma_2$ of~$\tsextic$ is a section of~$\Sigma_2$
disjoint from~$E$.
\pni
\endtheorem

\definition
Given a sextic $\sextic\subset\Cp2$ and a complete octuple
$\sQ_8\subset\DP(\sextic)$, the pair $(Y,\tsextic)$ equipped with the
projections $\Cp2\leftarrow Y\to\Sigma_2$ is called the \emph{\ACC-model}
of~$\sextic$ (defined by~$\sQ_8$). Here, $Y$ is a
relatively minimal rational Jacobian elliptic
surface and $\tsextic\subset Y$ is the bisection that projects
onto~$\sextic$.
\enddefinition

\subsection{The singular fibers}\label{s.fibers}

A complete octuple $\sQ_8\subset\DP(\sextic)$, see~\autoref{s.ACC},
is a union of maximal chains, one
chain (possibly empty) over each singular point~$P_i$ of~$\sextic$.
Given~$\sextic$, this octuple is determined by assigning the
\emph{height} $h_i:=\height P_i$
(the length of the corresponding chain) to each singular point~$P_i$. One has
$0\le h_i\le\frac12(p+1)$ if $P_i$ is of type~$\bA_p$ and $\sum h_i=8$.

\convention\label{conv.sextic}
Till the rest of this section, we fix a sextic $\sextic\subset\Cp2$
satisfying the following conditions:
\roster*
\item
$\sextic\subset\Cp2$ is a generic (see \autoref{def.generic})
irreducible non-special sextic,
\item
$\sextic$ has $\bA$-type singular points only and $\genus(\sextic)\le2$.
\endroster
Fix, further, a collection of heights $\{h_i\}$ of the singular points
of~$\sextic$ satisfying the conditions above. Hence,
we have also fixed a complete
octuple $\sQ_8\subset\DP(\sextic)$ and a pencil
$\pencil:=\pencil(\sQ_8)$ of cubics as in \autoref{s.ACC}, \ie, an \ACC-model
of~$\sextic$.
\endconvention

For a singular point~$P_i\in\sQ_8$,
we denote by $P_i\smax$ the topmost
(\ie, that of maximal level) element of~$\sQ_8$
that is infinitely near to~$P_i$.
If, in addition, $\height P_i\ge2$, then
$\fiber_i\in\pencil$ is the
member of the pencil singular at~$P_i$: such a
cubic obviously exists and, since a generic member of~$\pencil$ is
nonsingular, it is unique.

The following three statements are proved at the end of the section.

\theorem\label{th.fibers}
Each \emph{reducible} singular fiber of~$\t\pencil$
contains a \rom(unique\rom)
cubics~$\t\fiber_i$ corresponding to a singular point~$P_i$
of~$\sextic$ of height $h_i\ge2$.
\endtheorem

\addendum\label{ad.fibers}
Let $P_i$ be a singular point of~$\sextic$
and $h:=\height P_i\ge2$.
Then $\fiber_i$ is an irreducible nodal
\rom(possibly
cuspidal if $h\le3$\rom)
cubic.
The
corresponding fiber of $\t\pencil$ is
$\tF_i+\sum\tQ_j$,
the summation running over all points $Q_j\in\sQ_8$ such that
$Q_j\to P_i$ and $Q_j\ne P_i\smax$.
This fiber is of type~$\tA_{h-1}$,
possibly degenerating to~$\tA_{h}^*$
if $h\le3$.
\endaddendum

\addendum\label{ad.sections}
Let
$P_i$ be a singular point
of~$\sextic$
and $h:=\height P_i\ge1$.
Then $\tP_i\smax$ is a section of~$\t\pencil$ disjoint from~$\tQ_0$.
As a consequence,
the implicit basepoint~$Q_0$ of~$\pencil$ is a point of level zero.
\endaddendum

\proof[Proof of \autoref{th.fibers} and \autoref{ad.fibers}]
Let $\tF=m_1\tE_1+\ldots+m_r\tE_r$, $r\ge2$,
be a \emph{reducible} singular fiber.
Each component~$\tE_k$ of~$\tF$ is a nonsingular rational curve and one has
$\tE_k\cdot\tsextic\le2$, see \autoref{th.Artal}. Hence, $\tE_k$ lifts to a
nonsingular rational curve or a pair of such curves in the covering
$K3$-surface~$X$.

Assume that $\tE_k$ is not one of the exceptional divisors~$\tQ_j$, \ie,
$\tE_k$ projects to a curve $E_k\subset\Cp2$. Then $\deg E_k\le3$ and, due to
\autoref{th.rational}, the pull-back of~$\tE_k$ in~$X$ is irreducible.
Hence, $\tE_k\cdot\tsextic=2$ and, since also $\tF\cdot\tsextic=2$, we conclude that
$\tE_k$ is the \emph{only} component of~$\tF$ that does not contract to a
point in~$\Cp2$.
Thus, the image of~$\tF$ in~$\Cp2$
is either an irreducible rational cubic or a triple line.
In the former case, the singular point of the cubic should be resolved
in~$Y$; hence, this singular point is at one of those of~$\sextic$.
The latter case is easily ruled out as, due to \autoref{cor.generic}, the
line passes through two double points of~$\sextic$ only.

According to \autoref{th.rational}\iref{rational.3}, the cubic~$E_k$ passes
through seven double points $Q_1,\ldots,Q_7$ of~$\sextic$ and is singular at one of
these points, say, $P_i:=Q_1$. It is easy to see that $E_k$ belongs to~$\pencil$
if and only if all seven points are in~$\sQ_8$ and the eights element
of~$\sQ_8$ is infinitely near to~$P_i$.
Hence, the height of~$P_i$ is at least~$2$ and the corresponding singular fiber
of~$\pencil$ is as stated in \autoref{ad.sections}.
\endproof

\proof[Proof of \autoref{ad.sections}]
We have $K_Y\sim\tF$, where $\tF$ is a fiber of $\pencil$. Furthermore,
$(\tP_i\smax)^2=-2$ if $Q_0\to P_i\smax$ or
$(\tP_i\smax)^2=-1$ otherwise. By the adjunction
formula, in the former case $\tP_i\smax\cdot\tF=0$, \ie,
$\tP_i\smax$ is a component
of a (necessarily reducible) singular fiber of~$\t\pencil$. This possibility is
ruled out by \autoref{th.fibers} and \autoref{ad.fibers}.
In the latter case, $\tP_i\smax\cdot\tF=1$, \ie, $\tP_i\smax$ is a section.
Since $Q_0$ is not infinitely near to~$P_i\smax$, this section is disjoint from
$\tQ_0$.
\endproof

\subsection{The ramification locus \pdfstr{K in Sigma\sb2}{$\brl\subset\Sigma_2$}}\label{s.D}
Fix a sextic $\sextic\subset\Cp2$
and the other data as in \autoref{conv.sextic}
and
let $\brl\subset\Sigma_2$ be the trigonal part of the ramification locus of its
\ACC-model.

\corollary[of \autoref{th.fibers} and \autoref{ad.fibers}]\label{cor.D.sing}
The curve
$\brl\subset\Sigma_2$
has a type $\bA_{h_i-1}$ singular point $\barP_i$ for each
singular point~$P_i$ of~$\sextic$ of height $h_i\ge2$\rom;
this curve has no other singular points.
\done
\endcorollary

\proposition\label{prop.D.irreducible}
The
curve
$\brl\subset\Sigma_2$
is irreducible.
\endproposition

\proof
The curve~$\brl$ is reducible if and only if the Mordell--Weil group
$\MW(Y)$ has $2$-torsion,
see, \eg, \cite[Corollary 6.13 and Proposition 6.2]{degt:book}.
One has
\[*
\MW(Y)=H_2(Y)/\sset',\qquad
H_2(Y)=\Z[\t L]\oplus\bigoplus\Z[\tQ_i],
\]
where $\sset'\subset H_2(Y)$
is the sublattice
generated by the classes of the section and the components of all fibers
(see~\cite{Shioda:MW}), $L\subset\Cp2$ is a generic line,
and $Q_i\in\sQ_8^*$.
In view of \autoref{th.fibers} and
\autoref{ad.fibers},
we can decompose $\sset'=\sset''+\Z[\tF]$,
where $\sset''$ is generated by $[\tQ_0]$ and the
classes $[\tQ_i]$ of \emph{all but the topmost} elements $Q_i\in\sQ_8$
and $F$ is a generic fiber.
One has $[\tF]=3[\t L]-\sum(l_i+1)[\tQ_i]$,
where $Q_i\in\sQ_8^*$ and $l_i:=\QOPNAME{level}(Q_i)$.
Modulo~$\sset''$, the summation can be restricted to the topmost
elements $Q_i=P_j\smax$ only.
Then $l_i+1=h_j$ and,
since $\gcd\{h_j\}\mathrel|8$ is prime to~$3$,
the quotient $H_2(Y)/\sset'$ is torsion free.
\endproof

\corollary\label{cor.D.section}
Each singular point~$P_i$ of~$\sextic$
of height $h_i\ge1$
gives rise to a section $\sect_i:=\barP_i\smax$ of~$\Sigma_2$
\rom(\viz. the image of~$\tP_i\smax$\rom)
disjoint from~$E$.
This section is triple tangent to~$\brl$ \rom(if $h_i=1$\rom)
or
double tangent to~$\brl$ and
passing through $\barP_i$ \rom(if $h_i\ge2$\rom)\rom;
it does not pass through any other singular point of~$\brl$.
\endcorollary

In \autoref{cor.D.section}, we do not exclude the possibility that two or
three points of tangency of~$\barP_i\smax$ and~$\brl$ may collide.
Furthermore, these points of tangency may also collide with~$\barP_i$ (if
$h_i\ge2$).

\proof
In view of \autoref{prop.D.irreducible},
$\barP_i\smax$ is not a component of~$\brl$,
and the structure of the singular
fibers given by \autoref{ad.fibers} implies that $\barP_i\smax$ passes
through~$\barP_i$ (if $h_i\ge2$). Indeed,
in the notation of \autoref{ad.fibers},
the cubic~$\fiber_i$ does \emph{not} pass
through~$P_i\smax$, \ie, $\tF_i\cdot\tP_i\smax=0$.
On the other hand, it is $\tF_i$
that is the only component of the singular fiber that does not
contract in~$\Sigma_2$, see \autoref{ad.sections}.
Similarly,
$\barP_i\smax$ does not pass through any other singular point~$\barP_j$,
$P_j\ne P_i$, $h_j\ge2$, as the corresponding cubic~$\fiber_j$ \emph{does} pass
through~$P_i\smax$.
\endproof

\section{The Bertini involution}\label{S.Bertini}

The \ACC-construction can also be described as follows.
Let~$Y'$ be the plane blown up at the eight points $Q_1,\ldots,Q_8$. It is a
(nodal, in general) del Pezzo surface of degree~$1$.
The anti-bicanonical linear system defines a map
$\Gf\:Y'\to\Sigma_2'\subset\Cp3$, where $\Sigma_2'$ is a quadratic cone.
This map is of degree~$2$; it is ramified over the vertex
(the image of~$Q_0$) and a curve
$\brl'\subset\Sigma_2'$ cut off by a cubic surface disjoint from the vertex.
The image $\bsextic'\subset\Sigma_2'$ of~$\sextic$ is a plane section. The deck
translation of~$\Gf$ is called the \emph{Bertini involution} (defined by
the above pencil of cubics and its
distinguished basepoint~$Q_0$); it has one isolated fixed point, which is the
implicit basepoint~$Q_0$.
The objects appearing in the original construction, see \autoref{s.ACC}, are
obtained from those just described by blowing
this implicit point~$Q_0$ (or its image in~$\Sigma_2'$,
whichever is appropriate) up.

\subsection{Explicit equations\noaux{ (see~\cite{Moody})}}\label{s.Bertini}
In the exposition below, we try to keep the notation of~\cite{Moody}.
Consider the pencil defined by two plane cubics $\{w(x)=0\}$
and $\{w'(x)=0\}$, where
$x=(x_1:x_2:x_3)$, and assume that
$P_0(0:0:1)$, $P_1(0:1:0)$, and $P_2(1:0:0)$ are among its basepoints,
whence
\[*
w(x)=x_3^2(a_1x_1+a_2x_2)+
 x_3(b_1x_1^2+b_2x_1x_2+b_3x_2^2)+(c_1x_1^2x_2+c_2x_1x_2^2)
\]
and similar for~$w'$. The cubic passing through a point $y=(y_1:y_2:y_3)$ is
given by
\[*
\pW_3(x):=w(x)w'(y)-w'(x)w(y)=0.
\]
Clearly,
\[*
\pW_3(x)=x_3^2(A_1x_1+A_2x_2)+
x_3(B_1x_1^2+B_2x_1x_2+B_3x_2^2)+(C_1x_1^2x_2+C_2x_1x_2^2),
\]
where $A_i(y):=a_iw'(y)-a'_iw(y)$,
$B_i(y):=b_iw'(y)-b'_iw(y)$, and
$C_i(y):=c_iw'(y)-c'_iw(y)$.
Let $\kappa:=a_1b'_1 - a'_1b_1$ and consider the polynomials
\[*
\deMaple
\aligned
C_5(y)&:=A[2]*\[B[1] + \kappa*y[1]*y[3]^2]_{y[2]} +
  \[A[1] - \kappa*y[1]^2*y[3]]_{y[2]}*\[A[2]*y[3]+B[3]*y[2]]_{y[1]} +
  \kappa*B[3]*y[1]*y[3],\\
\phi_6(y)&:=A_1C_2+y[3]C_5,\\
\psi_6(y)&:=A_2C_1+y[3]C_5,\\
r_1'(y)&:=B_1A_2^2-B_2A_1A_2+B_3A_1^2.
\endaligned
\]
Here, following~\cite{Moody}, we use the notation $[e]_u$ to indicate that
$e$ has a common factor~$u$ and this factor has been removed.
In these notations, the Bertini involution is the birational map
$\Cp2\dashrightarrow\Cp2$, $y\mapsto z=(z_1:z_2:z_3)$, where
\[*
\deMaple
 z[1] = \phi_6*\[A[2]^2*\phi_6 + B[3]*r'[1]]_{y[1]},\quad
 z[2] = \psi_6*\[A[1]^2*\psi_6 + B[1]*r'[1]]_{y[2]},\quad
 z[3] = \psi_6*\phi_6*C_5.
\]
Apart from the basepoint~$P_0$, the fixed point locus is the order nine curve
$\rl\subset\Cp2$ given by the equation
\[*
\deMaple
\prl(y) := \psi_6*\[A[1]*y[3] + B[1]*y[1]]_{y[2]}
  -  \phi_6*\[A[2]*y[3] + B[3]*y[2]]_{y[1]} = 0.
\]
The sextics $\{\phi_6=0\}$ and $\{\psi_6=0\}$ play a special r\^{o}le: they are
the loci contracted by the Bertini involution to the basepoints~$P_1$
and~$P_2$, respectively.

\subsection{The map \pdfstr{Cp\sp2 --> Sigma\sb2}{$\Cp2\dashrightarrow\Sigma_2$}\noaux{ (see~\cite{degt:Bertini})}}\label{s.map}
The anti-bicanonical linear system $\ls|-2K_{Y'}|$ is generated by the strict
transforms of the sextics $\{\phi_6=0\}$,
$\{w^2=0\}$,
$\{ww'=0\}$, and $\{w\2=0\}$. Hence, in appropriate coordinates
$(\bz_0:\bz_1:\bz_2:\bz_3)$ in~$\Cp3$, the anti-bicanonical map $Y'\to\Cp3$,
regarded as a rational map $\Cp2\dashrightarrow\Cp3$, is given by
\[*
\bz_0=\phi_6(y),\quad
\bz_1=w^2(y),\quad
\bz_2=w(y)w'(y),\quad
\bz_3=w\2(y).
\]
Its image is the cone $\bz_1\bz_3=\bz_2^2$, which is further
rationally\mnote{`rationally'}
mapped to the
Hirzebruch surface $\Sigma_2$ \via\ $\bz\mapsto(\bx,\by)$, $\bx=\bz_1/\bz_2$,
$\by=\bz_0/\bz_2$. Here, $(\bx,\by)$ are \emph{affine} coordinates
in~$\Sigma_2$ such that the exceptional section~$E$ is the line
$\{\by=\infty\}$.
Summarizing, the composed rational map $\Cp2\dashrightarrow\Sigma_2$ defined
by a pair of cubics as in \autoref{s.Bertini} is $y\mapsto(\bx,\by)$,
\[*
\bx=w(y)/w'(y),\quad
\by=\phi_6(y)/w\2(y).
\]
Under this map, the pull-back of a proper section
\[
\{\by=\bpsextic_2(\bx)\},\quad
\bpsextic_2(\bx):=d_0+d_1\bx+d_2\bx^2
\label{eq.section}
\]
of~$\Sigma_2$ is the plane sextic~$\sextic$ given by the equation
\[
\phi_6(x)=d_0w\2(x)+d_1w'(x)w(x)+d_2w^2(x).
\label{eq.sextic}
\]
In particular, $\{\phi_6=0\}$ is the pull-back of the section $\{\by=0\}$.

The following conventions simplify the few further identities used in the
sequel.

\notation
Given a degree~$n$ monomial~$e$ in the coefficients
$a_1,\ldots,c_2$ of~$w$ and an integer $0\le m\le n$,
denote by $\(e)_m$ the sum of
$\binom{n}{m}$ monomials, each obtained from~$e$ by replacing $m$ of its~$n$
factors with their primed versions. For example,
\[*
\(a_1c_2)_1=a_1c_2'+a_1'c_2,\quad
\(b_2^2)_1=2b_2b_2',\quad
\(a_1b_1c_1)_2=a_1b_1'c_1'+a_1'b_1c_1'+a_1'b_1'c_1.
\]
This definition
extends to homogeneous polynomials by linearity.
\endnotation

\convention
Without further notice, we use same small letters to denote the coefficients of
a homogeneous bivariate polynomial:
$\pP_n(t_1,t_2)=\sum_{i=0}^np_it_1^it_2^{n-i}$ for a polynomial~$\pP_n$ of
degree~$n$.
With the common abuse of notation, we
freely treat homogeneous
bivariate polynomials as univariate ones: $\pP_n(\bx):=\pP_n(\bx,1)$.
This convention corresponds to the passage from homogeneous coordinates
$(t_1:t_2)$ to the affine coordinate $\bx:=t_1/t_2$ in the projective
line~$\Cp1$.
\endconvention

Since the strict transform of the sextic $\{\psi_6=0\}$ is also an
anti-bicanonical curve, there must be a relation $\psi_6=\phi_6+\pS_2(w,w')$
for a certain homogeneous polynomial~$\pS_2$ of degree two. Such a relation
indeed exists: one has
\[
s_0=a_2c_1-a_1c_2,\qquad
s_i=(-1)^i\(s_0)_i\quad\text{for $i=1,2$}.
\label{eq.pS}
\]
Hence, $\{\psi_6=0\}$ is the pull-back of the section $\{\by=-\pS_2(\bx)\}$
of~$\Sigma_2$.

Assume that the pencil has another basepoint
$P_3\notin(P_0P_1)$ of level zero and
normalize its
coordinates \via\ $u=(1:u_2:u_3)$. This point gives rise to another sextic
$\{\psi^u_6=0\}$, the one contracted to~$P_3$ by the Bertini
involution. Changing the coordinate triangle to $(P_0P_1P_3)$ and then
changing it back to $(P_0P_1P_2)$, one can easily see that
$\psi^u_6=\phi_6+\pS^u_2(w,w')$, where
\[
\gathered
\deMaple
s^u[0]=s[0]+(a[2]*c[2]*u_2+(a[2]*b[2]-a[1]*b[3])*u_3)+a[2]*b[3]*u_2*u_3+a[2]^2*u_3^2,\\
s^u_i=(-1)^i\(s^u_0)_i\quad\text{for $i=1,2$}.
\endgathered
\label{eq.pSu}
\]
As above, $\{\psi^u_6=0\}$ is the pull-back of the section
$\{\by=-\pS^u_2(\bx)\}$ of~$\Sigma_2$.

Finally, since $\rl\subset\Cp2$ is the pull-back of the ramification
locus $\brl\subset\Sigma_2$ (other than the exceptional section
$E\subset\Sigma_2$) and the latter is a proper trigonal curve,
there must be a relation
\[*
\prl^2=-4\phi_6^3+\phi_6^2\pP_2(w,w')+\phi_6\pQ_4(w,w')+\pR_3^2(w,w'),
\]
where $\pP_2$, $\pQ_4$, and $\pR_3$ are homogeneous polynomials of the
degrees indicated.
(The coefficient $(-4)$ is obtained by comparing the leading terms.)
The coefficients of these polynomials are
\[
\gathered
\deMaple
\aligned
 r[0] &= -a[1]*b[2]*c[2]+a[1]*b[3]*c[1]+a[2]*b[1]*c[2],\\
 q[0] &= 4*(a[1]*c[2] - b[1]*b[3])*s[0]+2*b[2]*r[0],\\
 p[0] &= b[2]^2-4*a[2]*c[1]-4*b[1]*b[3]+8*a[1]*c[2],
\endaligned\\
p_i=(-1)^i\{p_0\}_i,\quad
q_i=(-1)^i\{q_0\}_i,\quad
r_i=(-1)^i\{r_0\}_i\quad\text{for $i>0$}.
\endgathered
\label{eq.PQR}
\]
It follows that the defining equation of
the ramification locus $\brl\subset\Sigma_2$ is
\[
\bprl(\bx,\by):=-4\by^3+\by^2\pP_2(\bx)+\by\pQ_4(\bx)+\pR_3^2(\bx)=0.
\label{eq.K}
\]

\remark
Strictly speaking, most statements in
\autoref{s.Bertini} and \autoref{s.map}
hold only if
the pencil is sufficiently generic. Most important is the requirement that
the pencil should have no basepoints infinitely near to~$P_0$. Otherwise,
many expressions above acquire common factors; after the cancellation, the
Bertini involution degenerates to the so-called \emph{Geiser involution} and,
instead of a map $\Cp2\dashrightarrow\Sigma_2$, we obtain a generically
two-to-one map $\Cp2\dashrightarrow\Cp2$ ramified at a quartic curve (the
anti-canonical map of a nodal del Pezzo surface of degree~$2$).
See~\cite{Moody,degt:Bertini} for details.
Thus, in agreement with the \ACC-construction, we always assume that $P_0=Q_0$
is a simple basepoint of the pencil; the other basepoints may be multiple.
\endremark

\subsection{An implementation of the \ACC-construction}\label{s.implementation}
Together with Moody's formulas, the \ACC-construction gives us a relatively
simple way to obtain defining equations of sextics with large Milnor number.

Consider the pencil~$\pencil$ generated by a pair of cubics and, as in
\autoref{s.Bertini}, assume that it has at least three level zero
basepoints~$P_0$ (necessarily simple), $P_1$, $P_2$ at the coordinate
vertices and, possibly, some other basepoints $P_i$, $i\ge3$.
The pencil gives rise to a two-to-one rational map
$\Cp2\dashrightarrow\Sigma_2$, see \autoref{s.map}. We make use of the
following curves in~$\Sigma_2$:
\roster*
\item
the ramification locus $\brl=\{\bprl(\bx,\by)=0\}$, see~\eqref{eq.K};
\item
the section
$\sect_1:=\barP_1\smax=\{\by=0\}$, the image of $\{\phi_6=0\}$,
see \autoref{s.map}.
\item
the section
$\sect_2:=\barP_2\smax=\{\by=-\pS_2(\bx)\}$, the image of $\{\psi_6=0\}$,
see~\eqref{eq.pS};
\item
the sections
$\sect_i:=\barP_i\smax=\{\by=-\pS_2^u(\bx)\}$, $i\ge3$,
if present, see~\eqref{eq.pSu};
\item
the section $\bsextic=\{\by=\bpsextic_2(\bx)\}$, the image of
sextic~$\sextic$ to be constructed.
\endroster
Here, $\bpsextic_2$ is a degree~$2$ polynomial as in~\eqref{eq.section};
once found, it produces a sextic $\sextic\subset\Cp2$ given by
the defining equation~\eqref{eq.sextic}.

\remark
The pull-back of~$\sect_1$ in the elliptic surface~$Y$ splits into two
sections
interchanged by the deck translation. One of them projects
to~$\{\phi_6=0\}$, which is contracted to~$P_1$ by the Bertini
involution, see \autoref{s.Bertini}. Hence, the other is $\tP_1\smax$ and
$\sect_1$ is indeed~$\barP_1\smax$. The same argument applies to the other
sections~$\sect_i$.
\endremark

Let~$h_i$ be the multiplicity of the basepoint~$P_i$, $i\ge1$, \ie, the local
intersection index of the two cubics at this point. Then any sextic~$\sextic$
constructed as above is guaranteed to have a singular point adjacent
to~$\bA_{2h_i-1}$ at~$P_i$, and the construction is indeed the \ACC-model
of~$\sextic$ with $\height P_i=h_i$.
The further degeneration of~$\sextic$ depends on the position
of~$\bsextic$ with respect to~$\brl$ and~$\sect_i$, $i\ge1$. These
degenerations are discussed in details in \autoref{s.equations}.
Geometrically, the singularities of~$\sextic$ can be understood by running
the constraction backwards, \ie, considering the bisection
$\tsextic\subset Y$ and blowing it down to~$\Cp2$. Under the assumptions of
\autoref{conv.sextic}, the blow-down map $Y\to\Cp2$ is described by
\autoref{ad.fibers} and \autoref{cor.D.section}.
More precisely, for each $i\ge1$, we choose for~$\tP_i\smax$ one of the two
components of the strict transform of~$\sect_i$ in~$Y$, so that all
components chosen are pairwise disjoint.
(The existence of such a choice is guaranteed by the construction;
there are two coherent choices interchanged by the deck translation.)
If $h_i\ge2$, the section~$\sect_i$ passes through a singular point~$\barP_i$
of~$\brl$, which is in a certain singular fiber~$\bar\fiber_i$.
The corresponding
reducible
singular fibers of~$Y$ has several components:
the strict transform~$\tF_i$ of~$\bar\fiber_i$ and a number of other
components $\tQ_j$.
In these notations, the map $Y\to\Cp2$ is the blow-down of all
chosen sections~$\tP_i\smax$, followed by the consecutive blow-down of the
components~$\tQ_j\ne\tF_i$
of the reducible singular fibers, starting from the one
intersecting~$\tP_i\smax$. The components~$\tF_i$ left
uncontracted
project to the members
of the original pencil of cubics singular at~$P_i$, see \autoref{s.fibers}.

\figure[ht]
\centerline{\cpic{blow-down}}
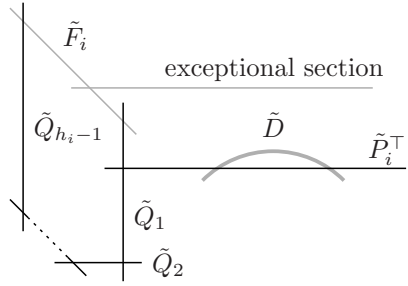
\caption{The divisors in~$Y$ blown down to $P_i\in\Cp2$}\label{fig.blow-down}
\endfigure

The divisors $\tP_i\smax,\tQ_j\subset Y$
blown down to a single singular point $P_i\in\Cp2$ are shown in
black in \autoref{fig.blow-down}, where the components $\tQ_j$ are numbered
in the order that they are contracted.

\remark\label{rem.parameters}
Since we are trying to find maximizing sextics, which are rare, we need
to consider families of pencils depending on parameters. It is easy to show
that the moduli space of the pencils that have basepoints of multiplicities
$h_0=1$, $h_1,\ldots,h_r$, $\sum h_i=9$,
has dimension~$r$. This agrees with the dimension
of the equisingular stratum of the moduli space of trigonal curves in~$\Sigma_2$:
assuming $\tA$ type singular fibers only, this dimension
equals $8-\mu(\brl)$.
\endremark

\remark\label{rem.forbidden}
Since we are interested in generic non-special irreducible sextics with double
singular points only, see \autoref{conv.sextic}, some values of the
parameters are forbidden, and we use this fact to simplify the equations.
Mostly, the following restrictions are used:
\roster
\item
The cubics $\{w=0\}$ and $\{w'=0\}$ are irreducible, see \autoref{th.fibers}
and \autoref{ad.fibers}. (In fact, \emph{all} members of the pencil must be
irreducible.)
\item
No three basepoints of level zero are collinear, see
\autoref{cor.generic}\iref{delPezzo.1}.
\item
The basepoint multiplicities are as stated, \ie, basepoints do not
collide.
\item
The section~$\bsextic$ is distinct from each~$\sect_i$ (as $\{\phi_6=0\}$ has
a triple point at~$P_1$, and similar for $\{\psi_6=0\}$ \etc.,
see~\cite{Moody}).
\item
More generally, $\bsextic$ does not pass through a singular point of~$\brl$,
as otherwise $\sextic$ would also have a triple singular point, \cf.
\autoref{fig.blow-down}.
\item
The section~$\bsextic$ is not tangent to
the ramification locus~$\brl$ within one of its reducible
singular fibers~$\bar\fiber_i$, see \autoref{cor.generic}\iref{delPezzo.3}.
\endroster
A number of other restrictions are ignored: we merely check the sextics
obtained and select those with the desired set of singularities.
\endremark

\section{The computation}\label{S.computation}

In this section, we outline some details of the computation.
Most polynomials
obtained
are too bulky to be reproduced here;
they can be downloaded from my web page~\cite{degt:equations},
in both human and machine readable form.
(I will extend this manuscript should there be any new development.)
Here, we provide information that is just enough to recover the
defining polynomials using
the formulas in \autoref{S.Bertini}.

We use the notation and setup introduced in \autoref{S.Bertini},
especially in \autoref{s.implementation}.

\subsection{Common equations}\label{s.equations}
We always assume that the multiplicities~$h_1$, $h_2$ of the
basepoints~$P_1$, $P_2$ are at least two and choose for
$\{w'=0\}$ and $\{w=0\}$ the unique members of the pencil that are singular
at~$P_1$ and~$P_2$, respectively.
Such pairs of cubics (with the necessary number of parameters)
are easily constructed by an appropriate triangular Cremona
transformation from appropriate pairs of conics.

Under the assumptions, the ramification locus $\brl\subset\Sigma_2$
has singular points $\bA_{h_1-1}$ and
$\bA_{h_2-1}$ over $\bx=\infty$ and $\bx=0$, respectively, and all sextics
obtained have singularities at least $\bA_{2h_i-1}$ at~$P_i$, $i=1,2$.
In the final equations, we change to affine coordinates $(x,y)$ in~$\Cp2$ so
that $P_1(0,\infty)$, $P_2(0,0)$, and the tangents to~$\sextic$ at these
points are the lines $\{x=\infty\}$ and $\{y=0\}$, respectively. This final
change of coordinates is indicated below for each pair of cubics.

The zero section $\sect_1$ intersects~$\brl$
at three double points, $\bx=\Gm,\Gn,\infty$, and in all
cases considered
(some of) the parameters present in
the equations are expressed rationally in terms of $\Gm,\Gn$.
For most equations, we use this re-parameterization.

The further degeneration of~$\sextic$ can be described using
\autoref{fig.blow-down} and the fact that each point of a $p$-fold, $p\ge2$,
intersection of~$\bsextic$ and~$\brl$ smooth for~$\brl$ gives rise to a
type~$\bA_{p-1}$ singular point of the strict transform $\tsextic\subset Y$.

If the section~$\bsextic$ is tangent to $\sect_1$,
\[*
\bpsextic_2(\bx)=a(\bx-\Gl)^2,\qquad
a\in\C^*,\quad
\Gl\in\C\sminus\{\Gm,\Gn\},
\]
the
$\bA_{2h_1-1}$ type singular point~$P_1$ of~$\sextic$
degenerates to $\bA_{2h_1}$.
If $\Gl=\Gm$, this point degenerates further to $\bA_{2h_1+1}$.
In this case, substituting $\by=a(\bx-\Gm)^2$ to the equation
$\bprl(\bx,\by)=0$
of the ramification locus, we obtain
\[
(\bx-\Gm)^2\pM_4(\bx-\Gm)=0,\quad
\pM_4(u):=m_0+m_1u+m_2u^2+m_3u^3+m_4u^4,
\label{eq.M}
\]
and the point~$P_1$ degenerates
to $\bA_{2h_1+1+k}$, $k\ge0$, if
\[
m_0=\ldots=m_{k-1}=0.
\label{eq.mk}
\]
The first equation $m_0=0$ is linear in~$a$; hence, $a$ can be
expressed rationally in terms of the other parameters and substituted to the
other equations.
Geometrically, this equation corresponds to the inflection tangency
of~$\bsextic$ and~$\brl$.

The degenerations of the other singular point~$P_2$
can be described similarly, by
analyzing the intersection of $\bsextic$ and~$\sect_2$.
Thus, the degeneration $\bA_{2h_2-1}\to\bA_{2h_2}$ is given by the
equation
\[
\QOPNAME{discriminant}(\bpsextic_2+\pS_2)=0.
\label{eq.second.point}
\]
Other singular points of~$\sextic$ are due to the extra tangency
of~$\bsextic$ and~$\brl$. Assuming that
$\bpsextic_2(\bx)=a(\bx-\Gm)^2$, the sextic has an extra singular point if
\[
\QOPNAME{discriminant}(\pM_{4-k})=0,\quad
\pM_{4-k}:=\pM_4/u^k.
\label{eq.discr}
\]
Note, though, that this discriminant may also vanish due to the
further
degeneration
$\bA_{2h_2-1}\to\bA_{2h_2+1}$ of~$P_2$ or another `fixed' singular point, if
present.
The sextic has an extra cusp (typically) if, in addition to~\eqref{eq.discr},
one has
\[
\QOPNAME{resultant}(\pM_{4-k},\pM''_{4-k})
=\QOPNAME{resultant}(\pM'_{4-k},\pM''_{4-k})=0.
\label{eq.cusp}
\]

\remark\label{rem.lost}
When simplifying equations and
their intermediate resultants, we routinely
disregard all factors that would result in forbidden (see
\autoref{rem.forbidden})
or
otherwise `unlikely'
values of
the parameters involved.
It is important to notice that,
since
the classification of sextics is known, we do not need to be too careful not to
lose a solution: it suffices to find the right number of distinct
curves realizing a given set of singularities.
The latter fact is given by \autoref{cor.all.curves} below.
\endremark

Typically, solutions to the equations appear in groups of Galois conjugate
ones: all unknowns are expressed as rational polynomials in a certain
algebraic number.
These groups are referred to as \emph{solution clusters}.

\subsection{The ramification locus \singset{2A3}}\label{s.2A3}
We start with a pair of cubics
\[*
\deMaple
\aligned
w &:= ((\Gb-\Ga)*x[1]+(\Gb-\Ga+\Ga*\Gb)*x[2]\strut)*x[3]^2\PLUS((\Gb-2*\Ga)*x[2]*x[1]-\Ga*x[2]^2)*x[3]-\Ga*x[1]*x[2]^2,\\
w' &:= ((\Ga-\Gb+\Ga*\Gb)*x[1]+(\Ga-\Gb)*x[2]\strut)*x[3]^2-(\Gb*x[1]^2-(\Ga-2*\Gb)*x[2]*x[1])*x[3]-\Gb*x[1]^2*x[2],
\endaligned
\]
where $\Ga,\Gb\in\C$, $\Ga,\Gb\ne0$, and $\Ga\ne\Gb$.
The change of coordinates for the final equations is
\[*
x_1 = 1-x,\qquad
x_2 = y-x,\qquad
x_3 = x.
\]
The re-parameterization in terms of $\Gm,\Gn$, see \autoref{s.equations}, is
as follows
\[*
\deMaple
\Ga = -\<(\Gm+1)(\Gn+1)//(\Gm+2)(\Gn+2)>,\qquad
\Gb = -\<(\Gm+1)(\Gn+1)//\Gm*\Gn+\Gm+\Gn+2>,
\label{eq.ab.2A3}
\]
and the equation $m_0=0$
(point~$P_1$ adjacent to \singset{A10}, see~\eqref{eq.mk}) yields
\[*
\deMaple
a=a[1] := \<(\Gm+\Gn+2)*(\Gm-\Gn)^2(\Gn+1)^2//4(\Gm+1)(\Gm+2)^2(\Gn+2)^4(\Gm*\Gn+\Gn+\Gm+2)>.
\]
All sextics below are obtained from sections~$\bsextic$ of the form
$\{\by=a_1(\bx-\Gm)^2\}$.

For the set of singularities \lineref{A12+A7}, the two additional equations
are $m_1=m_2=0$, see~\eqref{eq.mk}. They have two solutions
\[
\Gm = \frac1{13}(-4\pm6i),\quad
\Gn = 2\Gm-2,
\label{eq.A12+A7}
\]
producing two complex conjugate sextics.

For the set of singularities \lineref{A10+A7+A2}, the additional equations
are
\eqref{eq.discr} and~\eqref{eq.cusp}; their solutions are
\[
\Gm = 2\epsilon,\quad
\Gn = -10-12\epsilon,\qquad
\epsilon=\frac1{11}(-6\pm\sqrt3).
\label{eq.A10+A7+A2}
\]

Finally, assume $k=1$ in~\eqref{eq.mk} and consider
equations~\eqref{eq.second.point} and~\eqref{eq.discr}. They have three
solution clusters. One of them,
\[
\Gm = -5,\qquad
\Gn = -\frac{19}7,
\label{eq.D9+A10}
\]
results in a sextic with the set of singularities \singset{D9+A10}.
The two others are
\[
\Gm = \epsilon,\quad
\Gn = \frac15(-7+\epsilon),\qquad
\epsilon=\frac12(-11\pm3\sqrt5)
\label{eq.A10+A9}
\]
for the set of singularities \lineref{A10+A9} and
\[
\Gm = \frac\epsilon3,\quad
\Gn = \frac13(\epsilon^2+15\epsilon+13),\qquad
{\deZ _Z^3+17*_Z^2+51*_Z+43=0}.
\label{eq.A10+A8+A1}
\]
for the set of singularities \lineref{A10+A8+A1}.
In the former case, it is immediate that the line $y=0$ is not a component of
the curve; hence, the curves are irreducible.

\subsection{The ramification locus \singset{A5+A1}}\label{s.A5+A1}
We start with a pair of cubics
\[*
\deMaple
\aligned
w &:= (\Ga*x[1]-\Gb*x[2]\strut)*x[3]^2-(x[2]^2+(\Gb+1)*x[1]*x[2])*x[3]-x[1]*x[2]^2,\\
w' &:= ((\Gb+2*\Ga)*x[1]+\Ga*x[2]\strut)*x[3]^2\\
 \next+((2*\Ga+1)*x[1]*x[2]+(2*\Ga+\Gb+1)*x[1]^2)*x[3]+(\Ga+1)*x[1]^2*x[2],
\endaligned
\]
where $\Ga,\Gb\in\C$, $\Ga\ne0,-1$, and $\Ga+\Gb\ne0$.
The final change of variables is
\[*
\deMaple
x[1] = 1-x,\qquad
x[2] = y-\<(2*\Ga+\Gb+1)*x//\Ga+1>,\qquad
x[3] = x.
\]
The re-parameterization in terms of $\Gm,\Gn$, see \autoref{s.equations}, is
as follows
\[
{\deMaple
\Ga = 2-\Gm-\Gn,\qquad
\Gb = \Gm*\Gn-1,}
\label{eq.ab.A5+A1}
\]
and the equation $m_0=0$, see~\eqref{eq.mk}, result in
\[*
\deMaple
a = a_1 := -\<(\Gm+\Gn-2)*(\Gm-\Gn)^2//4(\Gm-1)^3(\Gn-2)>.
\]
The first three sextics are obtained from sections~$\bsextic$ of the form
$\{\by=a_1(\bx-\mu)^2\}$.

For the set of singularities \lineref{A16+A3}, the two additional equations
are $m_1=m_2=0$, see~\eqref{eq.mk}. They have two solutions
\[
\Gm = \epsilon,\quad
\Gn = \frac72-\epsilon,\qquad
\epsilon=\frac1{32}(59\pm3\sqrt{17}).
\label{eq.A16+A3}
\]
These curves were first studied in~\cite{Artal:trends}.

For the set of singularities \lineref{A15+A4}, the
equations
are $m_1=0$, see~\eqref{eq.mk}, and~\eqref{eq.second.point}.
They have two solutions
\[
\Gm = 1\pm3i,\quad
\Gn = \frac15(4-2\Gm).
\label{eq.A15+A4}
\]
These curves and their fundamental groups were studied
in~\cite{Artal:trends}.
It is easily seen that the conic maximally tangent to a curve at its
type~$\bA_{15}$ point is not a component. Hence, the curves are
irreducible.

Consider additional equations~\eqref{eq.second.point} and~\eqref{eq.discr}.
They have two solution clusters.
The first one,
\[
\Gm=4-3\Gn,\quad
\Gn=\frac16(5\pm i\sqrt{15}),
\label{eq.D5+A14}
\]
results in a sextic with the set of singularities \singset{D5+A14}.
The other one
\[
\gathered
\Gm = \epsilon,\quad
\Gn = \frac1{15}(9\epsilon^5-9\epsilon^4-72\epsilon^3+69\epsilon^2+172\epsilon-166),\\
\deZ 9*_Z^6-27*_Z^5-45*_Z^4+195*_Z^3-20*_Z^2-372*_Z+276=0
\endgathered
\label{eq.A14+A4+A1}
\]
gives us \lineref{A14+A4+A1}.
It is easily seen that $\Q(\epsilon)$ is a purely imaginary extension with
Galois group~$\DG{12}$.

Next two sextics are obtained from a section $\{\by=a_2(\bx-\Gm)^2\}$, where
\[*
\deMaple
a_2:=\<(\Gm+\Gn-2)*(\Gm^2*\Gn^2-6*\Gm*\Gn+4*\Gm+4*\Gn-3)//4(\Gm-1)^3>
\]
is found from equation~\eqref{eq.second.point}, which is linear in~$a$.
Adding equations~\eqref{eq.discr} and~\eqref{eq.cusp}, we obtain four
solution clusters.
One of them corresponds to a non-maximizing sextic, and another one,
\[
\Gm=5-3\Gn,\quad
\Gn=\frac16(7\pm i\sqrt3),
\label{eq.E6+A13}
\]
results in a sextic with the set of singularities \singset{E6+A13}.
The two others are
\[
\Gm = \epsilon,\quad
\Gn = \frac17(25-9\epsilon),\qquad
\epsilon=\frac12(7\pm\sqrt{21})
\label{eq.A13+A4+A2}
\]
for the set of singularities \lineref{A13+A4+A2} and
\[
\gathered
\Gm = \epsilon,\quad
\Gn = \frac18(9\epsilon^3-45\epsilon^2+73\epsilon-26),\\
\deZ 9*_Z^4-63*_Z^3+175*_Z^2-224*_Z+112=0
\endgathered
\label{eq.A13+A6}
\]
for the set of singularities \lineref{A13+A6}.
In the latter case, $\Q(\epsilon)$ is a purely imaginary extension with
Galois group~$\DG{8}$.

For the last set of singularities \lineref{A12+A6+A1},
substitution~\eqref{eq.ab.A5+A1} is not used.
A section $\bsextic=\{\by=a(\bx-\Gl)^2\}$ is tangent to~$\sect_2$
if and only if
\[*
\deMaple
a=a_2':=-\<\Ga*(4*\Ga+4*\Gb-\Gb^2)//4*(\Gl^2*\Ga+\Gl^2+\Gl*\Gb-2*\Gl+1)>,
\]
see~\eqref{eq.second.point},
the point of tangency being over
\[*
\deMaple
\bx=\Gl_2:=-\<\Gl*\Gb-2*\Gl+2//2*\Gl*\Ga+\Gb+2*\Gl-2>.
\]
Substitute $\by=a_2'(\bx-\Gl)^2$ to
$\bprl(\bx,\by)=0$
and expand the result as $\pN_6(\bx-\Gl_2)$, $\pN_6(u):=\sum_{i=0}^6n_iu^i$.
For the $\bA_6$ type point~$P_2$, we have $n_0=0$ (which
also implies $n_1=0$) and
$n_2=0$, and an extra $\bA_1$ type point results from the third equation
\[*
\QOPNAME{discriminant}(\pN_6/u^3)=0.
\]
This is a lengthy computation, and we use a certain `cheating', \cf.
\autoref{rem.lost}:
since the curves are expected to be defined over a cubic algebraic number
field, in all univariant resultants computed on the way we ignore all
irreducible factors of degree other than~$3$. At the end, we arrive at the
following solutions:
\[
\gathered
\Ga=4\epsilon,\quad
\Gb=-\frac1{27}(252\epsilon^2+468\epsilon+76),\quad
\Gl=-\frac1{27}(882\epsilon^2+1638\epsilon+329),\\
\deZ 441*_Z^3+315*_Z^2+79*_Z+7=0.
\endgathered
\label{eq.A12+A6+A1}
\]
Collecting all data, the sextic is the pull-back of the section
\[*
5103\by = (441\epsilon^2+204\epsilon+28)(27\bx+882\epsilon^2+1638\epsilon+329)^2.
\]

\subsection{The ramification locus \singset{A3+2A1}}\label{s.A3+2A1}
We start with a pair of cubics
\[*
\deMaple
\aligned
w &:= ((\Gb-\Ga)*x[1]+(\Ga*\Gb-2*\Ga+1\strut)*x[2])*x[3]^2\\
 \next+((\Ga*\Gb*\Gr-2*\Ga*\Gr+1)*x[1]*x[2]+(\Ga-\Ga*\Gr)*x[2]^2)*x[3]+(\Ga*\Gr-\Ga*\Gr^2)*x[1]*x[2]^2,\\
w' &:= ((\Ga*\Gb-2*\Gb+1)*x[1]+(\Ga-\Gb)*x[2]\strut)*x[3]^2\\
 \next+((\Ga*\Gb*\Gr-2*\Gb*\Gr+1)*x[1]^2+(-2*\Gb*\Gr+\Ga+\Ga*\Gr)*x[1]*x[2])*x[3]+(\Ga*\Gr-\Gb*\Gr^2)*x[1]^2*x[2],
\endaligned
\]
where $\Ga,\Gb,\Gr\in\C$, $\Ga\ne0,1$, $\Gb\ne0,1$, $\Gr\ne0,1$,
$\Ga\ne\Gb$, and $\Gb\Gr\ne1$.
The final change of variables is
\[*
\deMaple
x[1] = 1-\<x//\Gr>,\qquad
x[2] = y+\<(\Ga*\Gb*\Gr-2*\Gb*\Gr+1)*x//\Gr*(\Gb*\Gr-\Ga)>,\qquad
x[3] = x.
\]

This pencil of cubics has another $2$-fold basepoint $P_3(u)$,
\[*
\deMaple
u=(1:u_2:u_3):=(1:-\<\Gb*\Gr-1//\Gr-1>:-\<\Gb*\Gr-1//\Gb-1>),
\]
resulting in an $\bA_1$ type point of~$\brl$ over
$\bx=(\Gr-1)/(\Gb\Gr-1)$.
The corresponding
section is $\sect_3:=\barP_3\smax=\{\by+\pS_2^u(\bx)=0\}$,
see~\eqref{eq.pSu}.

The re-parameterization in terms of $\Gm,\Gn$, see \autoref{s.equations}, is
as follows
\[*
{\deMaple
\Gb = \<\Gm*\Gn*\Ga+2*\Gm*\Ga+2*\Gn*\Ga-\Gm-\Gn+3*\Ga-2//\Gm*\Gn+\Gm*\Ga+\Gn*\Ga+2*\Ga-1>,\qquad
\Gr = \<\Gm*\Gn+\Gm*\Ga+\Gn*\Ga+2*\Ga-1//(\Gm+\Ga)(\Gn+\Ga)>,}
\label{eq.ab.A3+2A1}
\]
and the equation $m_0=0$, see~\eqref{eq.mk}, results in
\[*
\deMaple
a = a_1 := \<\Ga*(\Ga-1)^4*(\Gm+\Gn+2)*(\Gm-\Gn)^2//
 4(\Gm+1)^2(\Gn+2)(\Gm+\Ga)^3(\Gn+\Ga)^2>.
\]

The first six sextics, with the sets of singularities adjacent to
\singset{A10+A4+A3}, are obtained from sections
$\bsextic=\{\by=a_1(\bx-\Gm)^2\}$, with $\Gb$, $\Gr$ as above and
\[*
\deMaple
\Ga=\<(\Gm*\Gn-1)^2*(\Gn+2)//4*\Gm*\Gn^2+10*\Gm*\Gn+\Gm^2+5*\Gn^2+8*\Gm+12*\Gn+8>
\]
found from~\eqref{eq.second.point}.
We need two more equations for the two parameters $\Gm,\Gn$ left.

For the set of singularities \lineref{A12+A4+A3}, the two
extra equations are
$m_1=m_2=0$, see~\eqref{eq.mk}. Their only solution is
\[
\Gm = -\frac{33}{13},\quad
\Gn = -\frac{29}{39}.
\label{eq.A12+A4+A3}
\]

For the set of singularities \lineref{A11+2A4},
we have
$m_1=0$, see~\eqref{eq.mk}, and
\[
\QOPNAME{discriminant}(\bpsextic_2+\pS^u_2)=0,
\label{eq.third.point}
\]
\cf.~\eqref{eq.second.point}.
These equations have four solutions
\[
\gathered
\Gm = -\frac1{11}(10\epsilon^3+70\epsilon^2+141\epsilon+109),\quad
\Gn = \epsilon,\\
\deZ 50*_Z^4+300*_Z^3+685*_Z^2+720*_Z+302=0.
\endgathered
\label{eq.A11+2A4}
\]
Analysing the discriminant, one can easily guess that the
splitting field of the above minimal polynomial
for~$\epsilon$ is $\Q(\sqrt2,i\sqrt{15})$.
An extra change of variables, making the two $\bA_4$ type points complex
conjugate, takes the four curves into two real ones, defined
and Galois conjugate over
$\Q(\sqrt2)$.

The pair of equations~\eqref{eq.discr} and~\eqref{eq.cusp}
has six solution clusters. One is
\[
\Gm = -\frac13,\quad
\Gn = \frac17(-14\pm i\sqrt7)
\label{eq.A10+A6+A3}
\]
for the set of singularities \lineref{A10+A6+A3}, and
another one is
\[
\Gm = \frac{17}{11},\quad
\Gn = \frac{41}{11}
\label{eq.A10+A4+A3+A2}
\]
for \lineref{A10+A4+A3+A2}.
Two solution clusters produce non-maximizing sextics, and the two others are
\[
\Gm = \frac16(-13\pm\sqrt{33}),\quad
\Gn = 3\Gm+2
\label{eq.E6+A10+A3}
\]
for the set of singularities \singset{E6+A10+A3} and
\[
\Gm = \epsilon,\quad
\Gn = \frac15(-\epsilon^2+4\epsilon+20),\qquad
{\deZ _Z^3+_Z^2-10_Z+10=0}
\label{eq.D5+A10+A4}
\]
for the set of singularities \singset{D5+A10+A4}.

The pair of equations~\eqref{eq.discr} and~\eqref{eq.third.point} has, among
others,
solutions
\[
\Gm = \frac13(-13+7\Gn),\quad
\Gn = \frac17(-11\pm3\sqrt{15}),
\label{eq.A10+A5+A4}
\]
resulting in the set of singularities \lineref{A10+A5+A4}, and
\[
\gathered
\deMaple
\Gm=-\<53025//51748>\Gn^5-\<14425//3044>\Gn^4-\<417315//51748>\Gn^3
 -\<440835//51748>\Gn^2-\<96864//12937>\Gn-\<59839//12937>,\\
\deMaple 875*\Gn^6+5375*\Gn^5+13375*\Gn^4+18025*\Gn^3+14770*\Gn^2+7180*\Gn+1592=0,
\endgathered
\label{eq.A10+2A4+A1}
\]
resulting in \lineref{A10+2A4+A1}.
The minimal polynomial for~$\Gn$ is reducible over $\Q(\Go)$,
$\Go:=i\sqrt{55}$. Furthermore, one can easily see that
\[
\Gn=\frac1{770}(21\Go\epsilon^2+87\Go\epsilon+385\epsilon+35\Go),\qquad
{\deZ 7*_Z^3+43*_Z^2+77*_Z+49=0},
\label{eq.A10+2A4+A1'}
\]
and another change of variables in~$\Cp2$,
making the two $\bA_4$ points
complex conjugate, converts the six sextics into three ones, defined
and Galois conjugate over $\Q(\epsilon)$.

The three other solutions to~\eqref{eq.discr} and~\eqref{eq.third.point}
provide an alternative representation of \singset{A10+A5+A4} and two
alternative representations for \singset{D5+A10+A4}.

\subsection{The ramification locus \singset{A3+2A1} (continued)}\label{s.A3+2A1.2}
For the remaining five sets of singularities, we start with the same pair of
cubics as in \autoref{s.A3+2A1} and observe that $\Gr$
can be expressed rationally in terms of the
$\bx$-coordinate~$\Gl$ of one of the points of tangency of~$\brl$
and~$\sect_3$:
\[*
\deMaple
\Gr = \<(1+\Gl*\Gb)*(\Gl*\Gb*\Ga+\Gl*\Ga-2*\Gl*\Gb+2*\Ga-\Gb-1)//
 \Gb*(\Gl+1)*(\Gl*\Gb*\Ga+\Gl*\Ga-2*\Gl*\Gb-\Gb*\Ga+3*\Ga-2)>.
\]
In all five cases, $P_3$ is adjacent to~$\bA_6$;
hence, $\bpsextic_2(\bx)=-\pS^u_2(\bx)+a(\bx-\Gl)^2$.
Then, substituting $\by=\bpsextic_2(\bx)$ to
$\bprl(\bx,\by)=0$,
we obtain
\[*
(\bx-\Gl)^2\pM_4(\bx-\Gl)=0
\]
for a certain polynomial~$\pM_4$ of degree four,
\cf.~\eqref{eq.M}, and the coefficient~$a$ is found from the equation
$m_0=0$. (The expression is too bulky to be reproduced here.)
In all five cases, we also have equation~\eqref{eq.second.point} making~$P_2$
adjacent to~$\bA_4$.

For the set of singularities \lineref{A8+A7+A4}, we have additional equations
\[
\QOPNAME{discriminant}(\bpsextic_2)=0
\label{eq.first.point}
\]
($P_1$ is adjacent to~$\bA_8$) and $m_1=0$ ($P_2$ is adjacent to~$\bA_7$).
The solutions are
\[
\Ga = \frac1{15}(27-14\epsilon),\quad
\Gb = -\frac1{45}(64+23\epsilon),\quad
\Gl = \frac1{37}(15+90\epsilon),\qquad
\epsilon=\pm i.
\label{eq.A8+A7+A4}
\]

In the other four cases, we use a `cheating' as above: since the curves are
expected to be defined over algebraic number fields of degree two or three,
we precompute univariate resultants and ignore their factors of degree greater than
four. (In the case \lineref{A7+2A6}, the presence of the two $\bA_6$ points
treated `asymmetrically' may
and does increase the field of definition.)

Equations~\eqref{eq.first.point} and~\eqref{eq.discr}, $k=1$, have four
solution clusters. One of them is
\[
{\deMaple\gathered
\Ga = \<1//2576595>(-820*\e9^2+559955*\e9+3862092),\quad
\Gb =\<\e9//9>,\\
\Gl=\<15//1098463796>(-44995*\e9^2+31556708*\e9-151837233),\\
\deZ 5*_Z^3-3495*_Z^2+8047*_Z-10925=0
\endgathered}\label{eq.A8+A6+A4+A1}
\]
for the set of singularities \lineref{A8+A6+A4+A1}, and another is
\[
{\deMaple\gathered
\Ga = \<1//226590>(-20121\e7^2+1110632\e7+22549),\quad
\Gb = 4*\e7,\\
\Gl= \<1//5395>(61959\e7^2-3470518\e7+41949),\\
\deZ 57*_Z^3-3196*_Z^2+221*_Z-7=0
\endgathered}\label{eq.A9+A6+A4}
\]
for the set of singularities \lineref{A9+A6+A4}.
The two others are
\[
\Ga = -\frac{13}7,\quad
\Gb = 91,\quad
\Gl = -\frac1{13}
\label{eq.D9+A6+A4}
\]
for the set of singularities \singset{D9+A6+A4} and
\[
\Ga = \frac1{21}(7\pm2i\sqrt7),\quad
\Gb = \frac1{66}(49\Ga-7),\quad
\Gl = -\frac34(21\Ga+11)
\label{eq.D5+A8+A6}
\]
for the set of singularities \singset{D5+A8+A6}.

Finally, consider~\eqref{eq.second.point} and the equations
\[*
\deMaple
3*m[2]*m[4] = m[3]^2,\qquad
3*m[1]*m[3] = m[2]^2,\qquad
9*m[1]*m[4] = m[2]*m[3].
\]
(This is a simplified version of~\eqref{eq.discr} and~\eqref{eq.cusp},
stating that a cubic polynomial is a perfect cube.)
They have three solution clusters:
\[
\Ga = \frac27(2\pm i\sqrt3),\quad
\Gb = \frac1{28}(19\Ga-6),\quad
\Gl = \frac14(7\Ga-22),
\label{eq.E6+A7+A6}
\]
resulting in the set of singularities \singset{E6+A7+A6},
\[
\Ga = \frac2{13}(40\Gb+9),\quad
\Gb = \frac1{50}(57\pm13\sqrt{21}),\quad
\Gl = \frac1{39}(25\Gb-22),
\label{eq.A7+A6+A4+A2}
\]
resulting in the set of singularities \lineref{A7+A6+A4+A2}, and
\[
\gathered
\Ga = \frac2{27}(7\Gb^3+42\Gb^2+66\Gb+20),\quad
\Gl=\frac1{54}(49\Gb^3+203\Gb^2+203\Gb-104),\\
49\Gb^4+245\Gb^3+357\Gb^2+56\Gb+22=0.
\endgathered\label{eq.A7+2A6}
\]
resulting in \lineref{A7+2A6}.
In this latter case, the minimal polynomial for~$\Gb$ becomes reducible over
$\Q(i\sqrt{7})$, and an extra change of variables converts the four curves
found into two complex conjugate curves defined over $\Q(i\sqrt{7})$.

\section{Proofs}\label{S.proofs}

In this concluding section, we outline the proofs of
the principal theorems stated in the introduction.
For \autoref{th.minimal}, we suggest two slightly different proofs. The one
in \autoref{proof.minimal} would work, with appropriate modifications and
computation, for any maximizing sextic, with or without triple singular
points. The other one, see \autoref{proof.alt}, is more limited, but it
reveals additional information (rather negative)
about the dessin of the cubic resolvent of a
sextic, see \autoref{rem.bad.dessin}.

\subsection{Proof of \autoref{th.minimal}}\label{proof.minimal}
The fields~$\Bbbk$ are described together with the equations of the curves,
and the computation of their Galois groups is straightforward.
For the minimality, we construct a projective invariant $J\in\C$
of (some) maximizing sextics, depending
rationally on the coefficients of their defining polynomials.

We use the following obvious observation: given a
polynomial
$d\in\C[x]$, the product of all linear factors of~$d$ of the same given
multiplicity is defined over the same field as~$d$.
Each curve in question has a distinguished singular point~$P_1$ of
type~$\bA_m$, $m\ge3$: for example, we can choose the point of the
maximal Milnor number.
Let $\psextic\in\Bbbk[x,y]$ be a defining polynomial of the curve
in some coordinate system.
Applying the above observation to the discriminant of~$\psextic$ with respect
to~$y$ (and then to the restriction to the corresponding fiber $x=\const$),
we conclude that the coordinates of~$P_1$ are in~$\Bbbk$. Hence, up to the
action of $\PGL(3,\Bbbk)$, we can assume that $P_1$ is $(\infty, 0)$ and that
the tangent to~$\sextic$ at~$P_1$ is the line $x=\infty$.

Let $d\in\Bbbk[x]$ be the discriminant of the new defining polynomial with
respect to~$y$. It has six distinct roots. (This is a common count for any
maximizing sextic with $\bA$-type singular points only, provided that all
its singular fibers are maximally generic.)
All curves considered have two to four singular points, each but~$P_1$
contributing a multiple root to~$d$.
Hence, $d$
has three to five simple roots; let
$\tilde d(x)=\sum d_ix^i$ be the product of the
corresponding linear factors. As explained above, $\tilde d\in\Bbbk[x]$,
and one can take for~$J$ any rational function of its coefficients~$d_i$
invariant under the action of the group
of affine linear transformations $x\mapsto ax+b$, $a\in\C^*$, $b\in\C$.
For our purposes, the following invariants are sufficient:
\roster*
\item
if $\deg\tilde d=3$, then
$J=J_3$
is the $j$-invariant of $\infty$ and the three roots of~$\tilde d$;
\item
if $\deg\tilde d=4$, then
$J=J_4$
is the $j$-invariant of the four roots of~$\tilde d$;
\item
if $\deg\tilde d=5$, then
$J=J_5:=\bigl(5d_3d_5-2d_4^2\bigr)^{10}\!/\bigl(5^{10}d_5^{12}\QOPNAME{discriminant}(\tilde d)\bigr)$.
\endroster
Explicitly, in the former two cases one has
\[*
\deMaple
J_3=-\<4*(3*d[1]*d[3]-d[2]^2)^3//27d[3]^2\Delta>,\qquad
J_4=\<4(d[2]^2+12*d[0]*d[4]-3*d[1]*d[3])^3//27\Delta>,
\]
where $\Delta:=\QOPNAME{discriminant}(\tilde d)$.
By construction, one has $J\in\Bbbk$ and $J$ does not depend on the choice of
coordinates or particular defining equation.

Now, on the case by case basis, one can check that, for each curve considered
in \autoref{S.computation}, the invariant~$J$ is well defined, \ie, the
singular fibers are maximally generic.
In most cases, the field~$\Bbbk$ obtained in the computation equals $\Q(J)$,
and this fact concludes the proof.
The exceptional cases are \lineref{A11+2A4},
\lineref{A10+2A4+A1}, and \lineref{A7+2A6}.
Each of these curves has a pair of isomorphic singular points, which are
treated asymmetrically by the construction, and the field obtained is twice
as large as predicted. In each case, an extra change of variables reduces the
field of definition to $\Q(J)$.
\qed

\corollary[of the proof]\label{cor.all.curves}
All sextics obtained in \autoref{S.computation} are pairwise distinct.
As a consequence,
all sextics listed in \autoref{tab.sextics} are present in
\autoref{S.computation}.
\endcorollary

\proof
Within each set of
singularities, the curves differ by the value of~$J$, regarded as a complex
number. (Since the curves are Galois conjugate, the values $J\in\Bbbk$ in the
abstract field of definition are equal.)\mnote{next remark removed as
obsolete}
\endproof


\subsection{An alternative proof \via\ \emph{dessins d'enfants}}\label{proof.alt}
In this section, we discuss another projective invariant $j_0:=j_0(\sextic)$
with the same property as above: $\Q(j_0)=\Bbbk$ is the minimal field of
definition of~$\sextic$. (This property is
easily verified by a direct case-by-case computation.)

As above, pick a distinguished singular point~$P_1$ of type~$\bA_m$,
$m\ge3$, for example, the one of the maximal Milnor number.
Consider the plane $\Cp2(P_1)$ blown up at~$P_1$:
it is a Hirzebruch surface~$\Sigma_1$,
and the strict transform of~$\sextic$ is a tetragonal curve intersecting the
exceptional section at a single point, which is a singular point of
type~$\bA_{m-2}$. Blowing this point up and blowing down the corresponding
fiber, we convert $\Cp2(P_1)$ to a Hirzebruch surface~$\Sigma_2$; the strict
transform of~$\sextic$ is a \emph{proper} (\ie, disjoint from the exceptional
section) tetragonal curve $\tsextic\subset\Sigma_2$.
This curve has a \emph{cubic resolvent} $C\subset\Sigma_4$:
if $\tsextic$ is given by a \emph{reduced} equation
\[*
\psextic(x,y):=y^4+p(x)y^2+q(x)y+r(x)=0,
\]
then $C$ is the proper trigonal curve given by
\[*
\pC(x,y):=y^3-2p(x)y^2+b_1(x)y+q(x)^2=0,\qquad
b_1:=p^2-4r.
\]
One can see that $C$ is equipped with a distinguished
section $L:=\{y=0\}$ that splits into two components in the covering elliptic
surface. Furthermore, $\tsextic$ is recovered from the pair $(C,L)$ uniquely
up to the transformation $(y,x)\mapsto(-y,x)$.

Associated to~$C$ is its \emph{functional $j$-invariant}
\[*
j(x)=\frac{4(p^2+12r)^3}{27\QOPNAME{discriminant}(\pC,y)};
\]
it is a rational map $\Cp1\to\C\cup\{\infty\}$.
The graph $j\1(\R\cup\{\infty\})\subset\Cp1$,
decorated
as
shown in \autoref{fig.j}, is called the \emph{dessin} of~$C$.
\figure[ht]
\centerline{\cpic{j}}
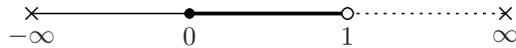
\caption{The decoration of the dessin}\label{fig.j}
\endfigure
This construction appears in a number of places; a detailed exposition
and further references can be
found in~\cite{degt:book}.
Typically, the \black- and \white-vertices of the dessin correspond to
\emph{nonsingular} fibers of~$C$ and have valency six and four, respectively,
whereas each \cross-vertex of valency~$2p$ corresponds to a singular fiber of
Kodaira's type~$\I_p$. The dessin may also have \emph{monochrome} vertices,
\viz.
the critical points of~$j$ with \emph{real} critical values other
that $0$, $1$, or~$\infty$.

It is easily seen that the total Milnor number of~$C$ is
$\Gm(C)=\Gm(\sextic)-2$ (assuming $\bA$ type singularities only).
Thus, if $\sextic$ is maximizing, $C$ is at most one unit short of being a
so-called \emph{maximal} trigonal curve, see~\cite{degt:book}. It follows
that~$j$ has at most one critical point with critical value
$j_0\ne0,1,\infty$; by definition, this critical value~$j_0$,
if defined, is the invariant being constructed.
\qed

\remark
The computation shows that the invariant~$j_0(\sextic)$ is well defined (and
has the property $\Q(j_0)=\Bbbk$) for all maximizing sextics with known
equations except \singset{(A17+A2)} (torus type) and \singset{A9+2A4+A2}
($\DG{10}$-special). In the two offending cases, the $\tA_2$ type singular
fiber of~$C$ degenerates to $\tA_2^*$ and $C$ is maximal.
\endremark

\remark\label{rem.bad.dessin}
The fact that $\Q(j_0)=\Bbbk$ proves also a certain negative result.
As explained above, typically,
the trigonal curve $C\subset\Sigma_4$ associated to a
maximizing sextic~$\sextic$ is almost maximal but not maximal.
On the other hand,
this curve is equipped with a distinguished section~$\sect$
splitting in the covering elliptic $K3$-surface, so that the latter has
maximal Picard rank, and this fact makes the pair $(C,\sect)$ rigid and
defined over an algebraic number field.
One might expect that the existence of such a section would manifest itself
in the combinatorial properties of the dessin of~$C$, \eg, in the presence of
a monochrome vertex, so that maximizing sextics with $\bA$ type singular
points only could also be studied in purely combinatorial terms.
However, this is not so: since $\Q(j_0)=\Bbbk$ is the minimal field of
definition, the only critical value~$j_0$ that could result in a monochrome
vertex is non-real whenever the chosen embedding $\Bbbk\into\C$ is non-real.
\endremark

\subsection{Proof of \autoref{th.pi1}}\label{proof.pi1}
All groups are computed as explained in~\cite{degt:tetra}, using real curves
and applying the Zariski--van Kampen theorem to real singular fibers only. In
the three exceptional cases, the presentations obtained are incomplete and
the computation is inconclusive.
Further details are found in~\cite{degt:equations}.

Two cases need special attention.
One is the set of singularities \lineref{A10+A7+A2}. This
set is realized by two real
sextics~$\sextic_1$, $\sextic_2$, one having one pair of
complex conjugate singular fibers, the other having two.
For the first curve, we have $\pi_1(\Cp2\sminus\sextic_1)=\CG6$. For the
other
one, the presentation is incomplete and the computation
only gives us an epimorphism
$G\onto\pi_1(\Cp2\sminus\sextic_2)$, where $G$ fits into a short exact
sequence
\[*
1\to\SL(2,\Bbb F_5)\to G\to\CG6\to1.
\]
(In other words, $[G,G]=\SL(2,\Bbb F_5)$, as found by \GAP.)\mnote{explained}
Hence, the group $\pi_1(\Cp2\sminus\sextic_2)$ is finite. On the other hand, since the two
curves are Galois conjugate, the profinite completions of their fundamental
groups are isomorphic, and we conclude that
$\pi_1(\Cp2\sminus\sextic_2)=\CG6$.

Alternatively, $\sextic_2$ can be projected from its $\bA_7$ type point. This
projection has only one pair of complex conjugate singular fibers, and we
obtain a complete presentation confirming that
$\pi_1(\Cp2\sminus\sextic_2)=\CG6$.

The other special case is the set of singularities \lineref{A10+A9} realized
by two real curves. For one of them, the presentation obtained using the
projection from the $\bA_{10}$ type point (as the equations suggest)
is incomplete. However, projecting
from the $\bA_9$ type point, we conclude that both groups are abelian.
\qed

\subsection{Proof of \autoref{th.Zariski}}\label{proof.Zariski}
We need to show that, within each pair $(\sextic_1,\sextic_2)$,
the spaces $\Cp2\sminus\sextic_i$, $i=1,2$, are not properly homotopy
equivalent.
(The sextics realizing
each of the sets of singularities
\singset{A18+A1} or \singset{A16+A2+A1} are also Galois conjugate;
this fact is proved
in~\cite{Artal:trends}.)

Let $\sextic\subset\Cp2$ be an irreducible sextic and consider the complement
$\Cp2\sminus\sextic$. Since
$H_1(\Cp2\sminus\sextic)=\CG6$, there is
a unique double covering
$X\open:=X\open_\sextic\to\Cp2\sminus\sextic$. It is an
oriented $4$-manifold;
hence, the group $H_2(X\open)$ is naturally equipped with the intersection
index form $H_2(X\open)\otimes H_2(X\open)\to\Z$.
The quotient $\bT_\sextic:=H_2(X\open)/\ker$
(where $\ker H_2=H_2^\perp$ stands for
the kernel of
the form) is a nondegenerate integral lattice.
Obviously, up to sign (a choice of orientation, \iq.
a generator for the group $H_c^4(X)=\Z$),
this lattice is preserved by proper homotopy equivalences
of the complement $\Cp2\sminus\sextic$.
In the case of maximizing sextics, the sign is determined by the requirement
that $\bT_\sextic$ should be positive definite.

Let $X\to\Cp2$ be the covering $K3$-surface of~$\sextic$, see \autoref{s.K3},
and let
$\sextic'\subset X$ be the preimage of~$\sextic$. Then
$X\open=X\sminus\sextic'$ and, by Poincar\'{e}--Lefschetz duality,
$H_2(X\open)=H^2(X,\sextic')$. From the exact sequence\mnote{ended with
$\ldots$}
\[*
H^1(X)=0\to H^1(\sextic')\to H^2(X,\sextic')\to H^2(X)\to H^2(\sextic')\to\ldots
\]
of pair $(X,\sextic')$ we conclude that
the invariant
$\bT_\sextic$ is
isomorphic to the orthogonal
complement of $H_2(\sextic')$ in $H_2(X)$.
Since $\sextic'$ contains all exceptional divisors and the divisorial
pull-back of $\sextic$ in~$X$ is equivalent to~$6h$,
the primitive hulls $H_2(D')\tilde\,$ and~$\tS$,
see \autoref{s.K3}, coincide; hence,
$\bT_\sextic=(\sset\oplus\Z h)^\perp$.
The latter orthogonal complement is often called the \emph{transcendental
lattice} of~$\sextic$; it plays an important r\^{o}le in the classification of
plane sextics.

Now, consulting Shimada's tables~\cite{Shimada:maximal}, one can see that, in
each pair as in the statement,
the two curves do differ by their transcendental lattices.
\qed

\subsection{The homotopy type of the complement}\label{s.homotopy}
In \autoref{proof.Zariski}, as well as in \autoref{proof.Zariski.weak} below,
we have to speak about \emph{proper}
homotopy equivalence only. The reason is the fact that
the definition of the intersection index form
uses Poincar\'{e} duality which, in the case of non-compact manifolds, involves
cohomology with compact supports. In general, this form does \emph{not}
need to be a homotopy invariant.
As an example, we show that the complements of most maximizing sextics
are homotopy equivalent.

Following~\cite{Dyer.Sieradski}, denote by~$P_m$ the
\emph{pseudo-projective plane}
of degree~$m$: this space is obtained by adjoining a $2$-cell~$e^2$
to the circle~$S^1$
\via\ a degree~$m$ map $\partial e^2\to S^1$.

\proposition\label{prop.homotopy}
Let $\sextic\subset\Cp2$ be a maximizing sextic with
$\pi_1(\Cp2\sminus\sextic)=\CG6$. Then there is a homotopy equivalence
$\Cp2\sminus\sextic\sim P_6\vee S^2$.
\endproposition

\proof
The complement $\Cp2\sminus\sextic$ is a Stein manifold; hence, it has
homotopy type of a {\sl CW}-complex of dimension~$2$.
Since the group $\pi_1(\Cp2\sminus\sextic)=\CG6$ is finite cyclic,
from~\cite{Dyer.Sieradski} one has
$\Cp2\sminus\sextic\sim P_6\vee S^2\vee\ldots\vee S^2$, where the number of
copies of~$S^2$ equals the Betti number $b_2(\Cp2\sminus\sextic)$.
Using Poincar\'{e}--Lefschetz duality, exact sequence of pair $(\Cp2,\sextic)$,
the fact that $\sextic$ is irreducible,
and the additivity of the topological Euler characteristic~$\chi$,
we obtain
\[*
b_2(\Cp2\sminus\sextic)
=b_1(\sextic)
=2-\chi(\sextic)
=2g-\Gm(\sextic)
=1,
\]
where $g=10$ is the genus of a nonsingular sextic.
\endproof

The proof of the following generalization is literally the same.

\proposition\label{prop.homotopy.general}
Let $\sextic\subset\Cp2$ be a plane curve of degree~$m$ with
$\pi_1(\Cp2\sminus\sextic)=\CG{m}$. Then
there is a homotopy equivalence
$\Cp2\sminus\sextic\sim P_m\vee S^2\vee\ldots\vee S^2$,
where the number of copies of the $2$-sphere~$S^2$ equals
$(m-1)(m-2)-\Gm(C)$.
\done
\endproposition

\autoref{prop.homotopy.general}
explains why $\pi_1$-equivalent Zariski pairs on
irreducible curves are so difficult to construct: if the
fundamental group is abelian,
the complements are not distinguished by the conventional homotopy
invariants.
In this respect, maximizing plane sextics are indeed very special, as their
transcendental lattices $\bT$ are positive definite and thus provide
additional invariants.
(The isomorphism class of an indefinite lattice is usually determined by its
signature and discriminant form, and these invariants can be computed in
terms of the combinatorial type.)
Another important class of curves with a similar property are the so-called
\emph{maximal trigonal curves} in Hirzebruch surfaces, see~\cite{degt:book}.

\remark
In a forthcoming paper, we will show that, with as few as about a dozen of
exceptions, the fundamental group of a non-special irreducible simple
sextic~$\sextic$
is~$\CG6$.
Hence, in most cases the homotopy type of
the complement $\Cp2\sminus\sextic$ is
completely determined
by~$\Gm(\sextic)$.
\endremark

\subsection{Proof of \autoref{th.Zariski.weak}}\label{proof.Zariski.weak}
%
Given a \emph{non-special} irreducible sextic~$\sextic$, in addition
to~$X\open$, see \autoref{proof.Zariski}, consider the double covering
$\barX\open\to\Cp2\sminus\Sing\sextic$ ramified at
$\sextic\sminus\Sing\sextic$.
Let $E\subset X$ be the union of the exceptional divisors in the covering
$K3$-surface, see \autoref{s.K3}.
Then we have Poincar\'{e}--Lefschetz duality $H_2(\barX\open)=H^2(X,E)$ and
exact sequence
\[*
H^1(E)=0\to H^2(X,E)\to H^2(X)\to H^2(E)=\sset.
\]
It follows that $H_2(\barX\open)=\sset^\perp$ is a non-degenerate lattice and
the inclusion $X\open\into\barX\open$ induces a primitive embedding
$\bT_\sextic\into\bS^\perp$, \cf. \autoref{proof.Zariski}.
Thus, the primitive lattice extension $\sset^\perp\supset\bT_\sextic$
(considered up to sign) is a
proper homotopy equivalence invariant of pairs~\eqref{eq.pairs}.
The orthogonal complement of~$\bT_\sextic$ in $\sset^\perp$ is $\Z h$, see
\autoref{s.K3}.

Now, using Nikulin's theory of discriminant forms, see~\cite{Nikulin:forms},
and the assumption that $\sset\oplus\Z h\subset\bL$ is a primitive
sublattice, see \autoref{th.special}, it is easy to show that the isomorphism
classes of primitive lattice extensions as above are
in a one-to-one correspondence with
the isomorphism classes of pairs $(\bT_\sextic,v\bmod2\bT_\sextic)$,
where $v\in\bT_\sextic$ is
a vector such that
$v\cdot\bT_\sextic\in2\Z$ and
$v^2=6\bmod8$. (The lattice $\sset^\perp$ is the index~$2$ extension
$(\bT_\sextic\oplus\Z h)+\Z v'$, where
$v':=\frac12(h+v)\in(\bT_\sextic\oplus\Z h)\otimes\Q$;
the arithmetical details are left to the reader.)

Using
Shimada's tables~\cite{Shimada:maximal} again,
we see that
$\bT_\sextic=\Z a\oplus\Z b$, where
\roster*
\item
$a^2=6$, $b^2=70$ for \lineref{A13+A4+A2},
\item
$a^2=22$, $b^2=30$ for \lineref{A10+A5+A4}, and
\item
$a^2=30$, $b^2=70$ for \singset{A6+A5+2A4}.
\endroster
In each case, there are two $\GROUP{O}(\bT_\sextic)$-orbits of classes
$(v\bmod2\bT_\sextic)$ as above, \viz. those of
$(a\bmod2\bT_\sextic)$ and $(b\bmod2\bT_\sextic)$.
These orbits distinguish the two sextics realizing~$\sset$.
\qed

\remark
According to~\cite{degt:JAG}, two simple sextics
$\sextic_1,\sextic_2\subset\Cp2$ have isomorphic oriented homological types
if and only if the pairs $(\Cp2,\sextic_i)$, $i=1,2$, are related by an
orientation preserving diffeomorphism subject to a certain regularity
condition at the singular points.
Diffeomorphisms of $4$-manifolds are a delicate subject, and the regularity
at the singular points is a subtle technical condition (roughly, it is
required that the structure of the exceptional divisors should be preserved).
The reader may observe that what is
essentially done in \autoref{proof.Zariski} and
\autoref{proof.Zariski.weak} is merely an attempt to extract
\emph{topological} invariants from the homological types.
\endremark

\let\.\DOTaccent
\def\cprime{$'$}
\bibliographystyle{amsplain}
\bibliography{degt}

\end{document}